\documentclass[12pt]{article}
\usepackage{subfiles}
\usepackage{float}
\usepackage{amsthm}
\usepackage{lineno}
\usepackage{cite}
\usepackage{amssymb,graphics,color,cite,amsmath}
\usepackage{subfig}
\usepackage{graphicx}
\usepackage{epsfig}
\usepackage{psfrag}
\usepackage{caption}
\usepackage{tabularx} 
\usepackage{booktabs} 
\usepackage{multirow} 
\usepackage[margin=0.75in]{geometry}
\usepackage{float}
\usepackage{afterpage}
\usepackage[linesnumbered,ruled,vlined]{algorithm2e}
\usepackage{amsfonts} 
\SetArgSty{textnormal}
\usepackage[noblocks]{authblk}
\usepackage[linesnumbered,ruled,vlined]{algorithm2e}
\SetKwRepeat{Do}{do}{while}%

\SetCommentSty{mycommfont}
\usepackage{graphics}
\usepackage{amsmath, amsthm, amssymb}
\usepackage{xifthen}
\usepackage{bbm, bm}
\usepackage[dvipsnames]{xcolor}
\usepackage{scalerel}
\newcommand{\bF}{{\bf F}}
\newcommand{\bB}{{\bf B}}

\newcommand{\bS}{{\bf S}}
\newcommand{\bR}{{\bf R}}

\newtheorem{thm}{Theorem}[section]

\theoremstyle{plainNoItalics}

\newtheorem{prop}[thm]{Proposition}

\title{Krylov-based Adaptive-Rank Implicit Time Integrators for Stiff Problems with Application to Nonlinear Fokker-Planck Kinetic Models
}

\author[a]{Hamad El Kahza}
\author[b]{William Taitano}
\author[a]{Jing-Mei Qiu}
\author[b]{Luis Chac\'on}

\affil[a]{\small{Department of Mathematical Sciences, University of Delaware, Newark, DE 19716}}
\affil[b]{\small{Theoretical Division, Los Alamos National Laboratory, Los Alamos, NM 87545}}

\begin{document}
\maketitle
\begin{abstract}
We propose a high order adaptive-rank implicit integrators for stiff  time-dependent PDEs, leveraging extended Krylov subspaces to efficiently and adaptively populate low-rank solution bases. This allows for the accurate representation of solutions with significantly reduced computational costs. We further introduce an efficient mechanism for residual evaluation and an adaptive rank-seeking strategy that optimizes low-rank settings based on a comparison between the residual size and the local truncation errors of the time-stepping discretization. We demonstrate our approach with the challenging Lenard-Bernstein Fokker-Planck (LBFP) nonlinear equation, which describes collisional processes in a fully ionized plasma. The preservation of {the equilibrium state} is achieved through the Chang-Cooper discretization, and strict conservation of mass, momentum and energy via a Locally Macroscopic Conservative (LoMaC) procedure. The development of implicit adaptive-rank integrators, demonstrated here up to third-order temporal accuracy via diagonally implicit Runge-Kutta schemes, showcases superior performance in terms of accuracy, computational efficiency, equilibrium preservation, and conservation of macroscopic moments. This study offers a starting point for developing scalable, efficient, and accurate methods for high-dimensional time-dependent problems.
\end{abstract}

{\bf Key words:} Adaptive-rank, extended Krylov based, implicit Runge-Kutta integrators, structure preserving, Lenard-Bernstein Fokker-Planck, local macroscopic conservation.

\section{Introduction}
\label{sec:intro}

Adaptive-rank representations of the solution to high-dimensional partial-differential equations (PDEs) have recently emerged as a viable solution strategy to ameliorate the so-called curse of dimensionality, whereby the computational complexity grows exponentially with the dimensionality of the problem. The approach postulates the solution as a sum of separable products of lower-dimensional functions along axes coordinates, with coefficients found dynamically in time along with the one-dimensional functions themselves. Such adaptive-rank representations are premised on the realization that, for many applications of interest, the rank of the solution $r$ (i.e., the number of coefficients needed for its accurate representation) is much smaller than the total number of degrees of freedom of the mesh resolution, resulting in significant computational savings. In fact, the storage complexity of low-rank representation of high dimensional functions scales as $O(dNr^2)$ in the tensor train format \cite{oseledets2011tensor}, with $d$ the dimensionality of the solution and $N$ the number of mesh points along a single direction. That is, the computational complexity scales linearly with dimensionality $d$ instead of exponentially, highlighting the curse-of-dimensionality-breaking potential of adaptive-rank methods.

However, realizing this potential requires a suitable strategy to adaptively track the lowest possible rank $r$ for an accurate representation of functions. 
Two main strategies have emerged in the literature for exploiting low-rank structures in time-dependent PDEs: the dynamic low rank (DLR) approach \cite{nonnenmacher2008dynamical, lubich2014projector, kieri2016discretized, ceruti2022unconventional, ceruti2022rank, dektor2021dynamic}; and the step-and-truncate (SAT) method, either for explicit schemes \cite{dektor2021rank, guo2022low, guo2023local, guo2022conservative} or for implicit schemes \cite{rodgers2023implicit, nakao2023reduced}.
Within the DLR framework,  differential equations 
are constructed to update the low-rank basis functions in each dimension through a projection (so-called \(K\) and \(L\) steps), after which a differential equation for the coefficient matrix (so-called
\(S\) step)
is formed through a projection onto the updated bases in all dimensions. These equations are then solved sequentially, either explicitly, or implicitly in the case of stiff PDEs.  A potential challenge associated with the DLR approach, however, in addition to the need to solve three systems of equations, is how to systematically formulate high-order time integrators able to couple  stiff and nonstiff terms. On the other hand, the SAT approach is built upon the traditional method-of-lines (MOL) spatial-temporal full discretization of PDEs and it can be naturally designed to be of high-order accuracy with a mixture of implicit and explicit treatments \cite{nakao2023reduced}.

Within the SAT approach, explicit adaptive low-rank methods evolve the numerical solution explicitly in time in a low-rank format, followed by an augmentation and a truncation step (using SVD) to discover the adaptive-rank representation of time-dependent PDE solutions. Its extension to implicit time discretizations, however, is far less straightforward, necessitating tailored strategies for the discovery of low-rank bases within a high-dimensional implicit framework. There has been significant interest in the literature on implicit low-rank methods: Venturi and collaborators \cite{rodgers2023implicit} recently proposed an efficient implicit tensor integrator via a direct application of the tensor train (TT) Generalized Minimal RESidual (GMRES) algorithm proposed in \cite{dolgov2013tt}. Nakao et al. \cite{nakao2023reduced} proposed a hybrid DLR-SAT approach in which a high-order implicit/explicit discretization of the matrix differential equations arising from multi-scale PDE discretizations are derived under the SAT framework, complemented with a predictor-corrector strategy to adapt the basis functions and their rank in the DLR spirit. Recently, Appelo et al. \cite{appelo2024robust} proposed an explicit prediction of the basis used to implicitly solve a reduced matrix system; an implicit solve is then applied to evolve the bases if the residual is not small enough.

In this study, we consider an implicit adaptive-rank integrator in which the rank discovery is performed via an extended Krylov subspace method.
Our motivating problem is the nonlinear Lenard-Bernstein Fokker-Planck (LBFP) collisional kinetic equation, given by \cite{filbet2022fokker}:
\begin{equation}
    \label{eq:FP}
    \frac{\partial f_{\alpha}}{\partial t} = 
    \sum^{N_s}_{\beta=1} \nu_{\alpha\beta} \nabla_v \cdot \left[
    D_{\alpha\beta} \nabla_v f_{\alpha} + \left(\vec{v} - \vec{u}_{\alpha \beta} \right) f_{\alpha}
    \right],
\end{equation}
which describes the collisional relaxation between $N_s$ plasma species. Here, $f_{\alpha}\left(\vec{v},t\right)$ is the particle distribution function for species $\alpha$, which is a function of the 
velocity space $\vec{v} \in {\mathbb R}^3$ and time $t\in{\mathbb R}_+$. The coefficient $\nu_{\alpha\beta}$ is the collision frequency between species $\alpha$ and $\beta$, $D_{\alpha\beta}$ is the diffusion coefficient (proportional to the temperatures of species $\alpha$ and $\beta$), $\vec{u}_{\alpha\beta}$ is the mixed drift velocity, and $n_{\alpha}$ is the number density, all of which are integrals of the particle distribution functions (to be defined precisely later, and which make the problem nonlinear). The LBFP equation is difficult to solve because it is nonlinear, stiff, it strictly conserves mass, momentum and energy, and it features a nontrivial null space (given by the Maxwellian distribution function with steady-state density, drift velocity, and temperature consistent with conservation of mass, momentum, and energy).

A simpler (linear) prototype problem for the above nonlinear model, also featuring strict conservation and a null space (with periodic boundary conditions), which we will consider to introduce concepts,
is the classical diffusion equation:
\begin{equation}
    \frac{\partial f }{ \partial t} =   \nabla_v \cdot (D \nabla f).
\label{eq: heat}
\end{equation}
Discretizing the above time-dependent PDEs on a tensor product velocity-space grid leads to the following {\em matrix} differential equation for the matrix element $F_{ij}$ [which is the discrete value of $f$ at the cell $(i,j)$]:
\begin{equation}
\frac{\partial \mathbf{F}}{\partial t} = D_1 \mathbf{F} + \mathbf{F} D_2^T := \mathbf{D}(\mathbf{F}). 
\label{Semi-Discrete}
\end{equation}
The matrix differential equation \eqref{Semi-Discrete} must be further discretized in time. Here, we consider an implicit temporal update to allow for a multiscale temporal integration of the diffusion equation. The simplest example of an implicit scheme is the first-order backward Euler discretization of \eqref{Semi-Discrete} (although we will consider up to third-order Diagonally Implicit Runge-Kutta (DIRK) schemes later in this study), leading to the following Sylvester equation for $\mathbf{F}^{(n+1)}$, where the superscript $(n+1)$ denotes the time level of numerical solution with time stepping size $\Delta t$: 
\begin{align}
    \left(\frac{1}{2}I - \Delta t D_1\right)\mathbf{F}^{(n+1)} + \mathbf{F}^{(n+1)}\left(\frac{1}{2}I - \Delta t D_2^T\right) = \mathbf{F}^{(n)}.
    \label{Sylv_disc}
\end{align}

In this study, we uncover the low-rank structure of the solution to Eq. \eqref{Sylv_disc} dynamically in time using extended Krylov-subspace methods. Krylov-subspace methods are a class of iterative techniques for solving large linear systems of equations. They have proven to be powerful, especially when dealing with sparse matrices. The foundational work by Saad in 1989 \cite{saad1989numerical} laid the groundwork for using Krylov spaces to solve the Sylvester equation, which is a linear matrix equation that could arise in implicit MOL discretization of PDEs such as Eq. \eqref{Sylv_disc}. Saad's method looks for solutions with basis from each dimension constructed from Krylov subspaces, followed by a Galerkin projection to update the coefficient matrix from a reduced Sylvester equation. Built upon Saad's methodology, Simoncini in 2007 \cite{simoncini2007new} introduced the use of extended Krylov subspaces to enhance the convergence of solutions to matrix equations. This extension provides a more robust framework by combining both the Krylov space of a matrix and its inverse, thus generating a richer subspace that often leads to faster convergence for many problems. 
In our low-rank context, the extended Krylov method is rendered competitive because matrices are sparse, and inverted only in a single dimension at a time, which scales linearly with the number of dimensions $d$ and the one-dimensional mesh size $N$. 

We propose a criterion to stop the rank augmentation process adaptively by comparing the residual magnitude to the local truncation error (LTE) of the DIRK scheme, striking a balance between computational efficiency and accuracy. We estimate the residual norm while avoiding forming the full basis explicitly using an efficient and adaptive approach proposed by Shankar \cite{shank2013krylov}. These innovative developments coalesce into a framework capable of finding low-rank, high-order implicit solutions for stiff time-dependent parabolic PDEs with super-optimal scaling and at a significantly reduced computational cost vs. the full-rank algorithm. 

This work presents several key contributions to the field:
\begin{itemize}
\item Under the framework of developing efficient implicit adaptive-rank integrators, we leverage the extended Krylov subspaces to populate low-rank solution basis,
which provides a fertile ground for seeking implicit low-rank solutions for time-dependent problems. 
We then employ an efficient mechanism for residual evaluation and introduce an adaptive-rank-seeking approach. This strategy compares the residual size with the local truncation errors inherent in the time-stepping discretization, optimizing the low-rank settings to achieve desired accuracy with reduced computational complexity. The proposed algorithm is super-optimal in that it scales linearly with respect to resolution in a single dimensiona $N$ and the dimensionality $d$ (instead of $N^d$).
\item 
We apply the proposed time-dependent low-rank algorithm to the nonlinear Fokker-Planck collisional equation. We linearize it by separately evolving equations for mass, momentum and energy moments. We discretize it with the Chang-Cooper discretization \cite{chang1970practical}, which preserves the Maxwellian equilibrium analytically. For strict conservation of collisional invariants, we apply a Locally Macroscropic Conservative (LoMaC) procedure~\cite{guo2023local} to project the low-rank kinetic solution to a reference manifold defined by the macroscopic moments, conserving the mass, momentum, and energy up to machine precision.

\item 
We develop implicit 
adaptive-rank integrators 
via diagonally implicit Runge-Kutta schemes, which could be designed for arbitrarily high temporal order. Their performance in temporal accuracy (up to third order in this study), computational efficiency, equilibrium preservation, and macroscopic conservation are numerically demonstrated.  
\end{itemize}
The remainder of this paper is structured as follows. In Section \ref{sec:extended-krylov}, we introduce the proposed extended-Krylov-based low-rank implicit solver. In Section \ref{sec:fp} we discuss the multispecies nonlinear Fokker-Planck equation of interest, along with its conservation properties. We present numerical results demonstrating the properties of the algorithm in Sec. \ref{sec:num-exp}, and we conclude in Sec. \ref{sec:conclusions}.

%
%
\section{Extended Krylov adaptive-rank implicit integrators for stiff problems}
\label{sec:extended-krylov}
In this section, we discuss Krylov-based implicit low-rank algorithms for the classical heat equation \eqref{eq: heat} as a prototype problem for the general nonlinear Fokker-Planck model, which will be discussed later in Sec. \ref{sec:fp}. In Section~\ref{sec: 2.1}, we set up the classical MOL discretizations for the heat equation on a tensor product of 1D grids. In Section~\ref{sec: 2.2}, we propose a Krylov-based low-rank solver for the linear matrix equation, e.g. \eqref{Sylv_disc}. In Section \ref{sec: 2.3}, we analyze the computational complexity of the algorithm, and in Section~\ref{sec: 2.4}, we extend the proposed algorithm to high-order DIRK time-integration methods. 

%
%

\subsection{Adaptive-rank implicit integrators for the heat equation: basic setup}
\label{sec: 2.1}

We consider a two-dimensional tensor product grid and set the number of grid points in the $v_1$ and $v_2$ direction as ${N}_1$ and ${N}_2$, respectively.
We assume a low-rank approximation to the initial condition $\bF_0 \in \mathbb{R}^{{N}_1 \times {N}_2}$ {at time $t^{(0)}$, which we evolve to time $t^{(1)}=t^{(0)}+\Delta t$,} 
\begin{equation}
\bF_0 = U_0S_0V_0^T,
\end{equation}
where $U_0 \in \mathbb{R}^{{N}_1 \times r}$ and $V_0 \in \mathbb{R}^{{N}_2 \times r}$ with their orthonormal columns representing bases in the respective dimensions. $S_0 \in \mathbb{R}^{r \times r}$ is a diagonal matrix with decreasing singular values, which are coefficients for the outer product of basis functions. 

To discretize Eq. \eqref{eq: heat}, we follow the classical method-of-lines approach by first discretizing in velocity space and then in time. For the velocity-space discretization, we consider second-order finite differences, which generate tridiagonal matrices per dimension, optimizing computational efficiency by exploiting their sparse structure.
For the temporal discretization, we consider implicit schemes to address the numerical stiffness of the diffusion operator. The classical first-order scheme is the backward Euler method, the discretization of which leads to a linear matrix equation of the Sylvester type, e.g. Eq. \eqref{Sylv_disc}. High-order extensions, such as DIRK methods are possible and will be discussed in Section \ref{sec: 2.4}. Below, we consider a linear Sylvester equation of the form: 
\begin{equation}
      A_1 {\bF} + {\bF} A_2^{T} = {\bf B}
        \label{eq: Syl}
\end{equation}
where $A_1$ and $A_2$ are sparse discretization matrices applied to their respective dimensions, $\bF$ is the implicit update we are seeking, and $\bB$ depends on solutions at previous time steps. For example, for BE:
\begin{equation}
A_1=\frac{1}{2}I-\Delta t D_1, \  A_2=\frac{1}{2}I-\Delta t D_2, 
\label{diff_matrices}
\end{equation}
$\bF = \bF^{(n+1)}$, and $\bB=\bF^{(n)} =  U_0 S_0 V_0^T$ in \eqref{Sylv_disc}. When high-order DIRK methods are considered, $\bB$ can be explicitly evaluated from low-rank solutions at previous RK stages. 
Efficient solvers for the Sylvester equation \eqref{eq: Syl} have been developed in the past few decades \cite{saad1989numerical, simoncini2007new, simoncini2016computational}. Below, we leverage these developments and propose an efficient adaptive-rank extended Krylov-based implicit integrators for stiff time-dependent problems. 

%
%
\subsection{Extended Krylov methods for the Sylvester equation}
\label{sec: 2.2}

Reference \cite{simoncini2016computational} provides an extensive review on solving the Sylvester equation \eqref{eq: Syl}. The equation admits a solution if and only if the spectrum of $A_1$ and $-A_2$ are well separated. In particular, 
\begin{equation}
    \|\mathbf{F}\|_F \leq \frac{\|\bB\|_F}{\text{sep}({A_1}, -{A_2})}, 
\end{equation}
with $\text{sep}(A_1, -A_2) = \min_{\|P\|_p=1} \|A_1 P + P A_2\|_p$.
Equation \eqref{Sylv_disc} provides an example of the Sylvester equation arise from numerical discretization of diffusion equations.
The solution of \eqref{eq: Syl} admits the following form
    \begin{equation}
        \bF = -\int_0^{\infty} e^{A_1 \tau}\bB  e^{A_2^T \tau} d\tau \overset{\bB \doteq U_0 S_0 V_0^T}{\hstretch{3}{=}} 
-\int_0^{\infty} (e^{A_1 \tau} U_0) S_0   {(e^{A_2 \tau} V_0)}^T d\tau.
        \label{eq: exp_formu}
    \end{equation}
    Here $e^{A_1 \tau}U_0$ and $e^{A_2 \tau}V_0$ can be computed as
    \[
e^{A \tau}U {\approx} U + \tau AU + \frac{\tau^2}{2!} A^2U + \ldots +  \frac{\tau^{m-1}}{{m-1}!} A^{m-1}U = \kappa_m(A,U)Z_1^T.
\]
Here, $m$ is an integer such that $m \leq N$. This suggests the Krylov subspaces $\kappa_m(A_1,U_0)$ and 
$ \kappa_m(A_2,V_0)$ as basis candidates for the matrix solution $\bF$ in their respective dimension, with:
\begin{equation}
    \kappa_m(A,U) = [U, A U, A^2 U, \ldots, A^{m-1} U]. 
\end{equation}
Following this insight, a standard Krylov iterative method was proposed for linear matrix equation in Ref. \cite{saad1989numerical}. The method augments Krylov subspaces in an iterative fashion, until the residual norm is below a prescribed threshold tolerance. However, such method did not find practical success due to its slow convergence. Recent developments have led to the emergence of an extended Krylov iterative method with better convergence properties \cite{druskin1998extended, simoncini2007new}. 
The idea of extending Krylov subspaces stems from the equivalent formulation of the problem \eqref{eq: Syl} by multiplication of $A_1^{-1}$ and $A_2^{-1}$ from left and right respectively, assuming that $A_1$ and $A_2$ are invertible. That is, 
\[
\bF A_2^{-1}  + A_1^{-1} \bF  = (A_1^{-1}U_0) S_0 (A_2^{-1} V_0)^{T}. 
\]
Solving the above equation by the standard Krylov iterative method requires searching for an approximation in inverted Krylov subspaces given by: 
\[
\kappa_m(A_1^{-1},A_1^{-1}U_0) = [A_1^{-1}U_0,\ldots, A_1^{-m+1} U_0],
\qquad
\kappa_m(A_2^{-1},A_2^{-1}V_0) = [A_2^{-1}V_0,\ldots,A_2^{-m+1} V_0].
\]
The extended Krylov iterative method searches the solution from the following richer subspaces in their respective dimension $\kappa_m(A_1,A_1^{-1},U_0)$ and $\kappa_m(A_2,A_2^{-1},V_0)$ with:
\begin{equation}
\kappa_m(A,A^{-1},U) = [U, A U,A^{-1}U,\ldots, A^{m-1} U,A^{-m+1} U].
\label{eq: extend_K}
\end{equation}
Theoretical results concerning the convergence rate of extended Krylov subspaces to the action of matrix functions can be found in \cite{knizhnerman2010new}, along with comprehensive references. Motivated by these developments, we propose to explore the use of extended Krylov-based low-rank implicit integrators for stiff PDEs. In particular, we propose to use extended Krylov methods to adapt the implicit solution rank in time. These methods are practical in our context because the matrix inverses required are tridiagonal, and therefore fast to compute. The proposed  extended-Krylov implicit algorithm comprises the following steps:

\begin{enumerate}
    %
    %
    \item [Step K1.] {\em Construction of dimension-wise Krylov subspaces.}
We consider the following dimension-wise Krylov subspaces 
\begin{equation}
\kappa_m(A_1,A_1^{-1},U_0), \qquad \kappa_m(A_2,A_2^{-1},V_0)
\end{equation} 
for the respective dimension, with 
\begin{equation}
\label{eq: K_N1}
\kappa_m(A,A^{-1},U) = [U, A U,A^{-1}U,\ldots, A^m U,A^{-{m+1}} U].    
\end{equation}
Upon building the extended Krylov subspaces for each dimension, a reduced QR decomposition can be performed 
\begin{equation}
    \kappa_m(A_1, A_1^{-1}, U) = U^{(m)} R_U, \quad \kappa_m(A_2, A_2^{-1}, V) = V^{(m)} R_V,
\end{equation}
where {$U^{(m)}$ and $V^{(m)}$} form sets of orthonormal basis for Krylov subspaces, and $R_U$ and $R_V$ are upper triangular matrices defining the mapping from Krylov subspaces to their respective orthonormal basis.

    %
    %

\item [Step K2.] {\em Projection method for a reduced Sylvester equation.} Now that we have obtained the orthonormal Krylov bases $U^{(m)}$ and $V^{(m)}$, we seek a solution in the form of $\bF^{(m)} =  U^{(m)} S^{(m)} {V^{(m)}}^T$. In the following, we skip the upper script ${(m)}$ 
and let $\bF =  U_1 S_1 V_1^T$ be the evolved solution at time $t^{(1)}$, for notation simplicity.  We derive a reduced Sylvester equation for $S_1$ via a Galerkin projection of the residual,
\begin{equation}
    \bR= A_1 {\bF} + {\bF} {A_2^T}-U_0 S_0 V_0^T, \qquad \mbox{with} \quad \bF =  U_1 S_1 V_1^T.
    \label{residual}
\end{equation}
In particular, with Galerkin projection $(U_1)^T \bR V_1 = 0$, we have 
\begin{align}
       & \tilde{A_1} S_1 + S_1 \tilde{A_2}^T= \tilde{B}_1      , 
    \label{Galerkin}
\end{align}
where 
\begin{equation}
    \tilde{A_1} = U_1 ^{T} A_1 U_1 = \frac12 I - \Delta t U_1 ^{T} D_1 U_1 \doteq \frac12 I - \Delta t \tilde{D}_1,
    \label{eq: tA1}
\end{equation}
similarly for $\tilde{A_2}$, and 
\begin{equation}
    \tilde{B}_1 \doteq (U_1^T U_0) S_0 (V_1^T V_0)^T.
    \label{eq: bB1}
\end{equation}
The above reduced Sylvester equation is of the size of the Krylov subspaces, from which $S_1$ is obtained using a direct solver.

    %
    %

\item [Step K3.] {\em Efficient evaluation of the residual norm.}
In an iterative process, the Krylov iteration incrementally grows the size of Krylov subspaces \eqref{eq: K_N1}, based on evaluating the norm of the residual matrix $\lVert \bR\rVert_F$ in an efficient recursive fashion:
\begin{align}
\lVert \bR\rVert_F&= \lVert  {A_1}  {\bF}_{1} +   {\bF}_{1}  {A_2}^T - U_0 S_0 V_0^T\rVert_F \nonumber\\ 
\nonumber &=
\left\| 
\begin{bmatrix}
      {U}_{1} &  {A_1} 
  {U}_{1}
\end{bmatrix}
\begin{bmatrix}
    -({U}_1)^T  {U_0 S_0 V_0^T}  {V}_1 &  {{S_{1}}}  \\
     {{S_{1}}} & {\mathbf{0}} \\
\end{bmatrix}
\begin{bmatrix}
  {V}_{1}&  {A_2}  {V}_{1}
\end{bmatrix}^T
\right\|\\ \nonumber
=&\left\| 
 {Q}_U  {R}_U
\begin{bmatrix}
    - \tilde{B}_1 &  {{S_{1}}}  \\
     {{S_{1}}} & {\mathbf{0}} \\
\end{bmatrix}
 {R}_V^T  {Q}_V^T
\right\|\\ 
=& \left\| 
 {R}_U
\begin{bmatrix}
     -\tilde{B}_1 &  {{S_{1}}}  \\
     {{S_{1}}} & {\mathbf{0}} \\
\end{bmatrix}
 {R}_V^T \right\|, 
 \label{Residual}
\end{align}
where the second equality above results from $U_1 {U}_1^T U_0 S_0 V_0^T {V}_1 V_1^T = U_0 S_0 V_0^T$, where we have used that $U_1 {U}_1^T$ is a projection onto its column space, which includes $U_0$ due to the Krylov subspace construction, and similarly with $V_1 V_1^T$ and $V_0$. The factors
$Q_U R_U$ and $Q_VR_V$ are obtained using the reduced QR decomposition of matrices $[U_1, A_1 U_1]$ and $[V_1, A_2 V_1]$, respectively. As we argue in Sec. \ref{sec: 2.3} below, the computational complexity of this step does not scale with the problem size in a single dimension $N$, only with the cube of the low-rank dimension $r$, and is therefore inexpensive when $r$ remains small.

    %
    %

\item [Step K4.] {\em Adaptive augmentation of Krylov subspaces.} 
Once the residual norm is computed in \eqref{Residual}, we compare it with an estimate of the temporal discretization error to adaptively determine the size of the Krylov subspaces. In particular, the tolerance $\epsilon_{tol}$ for residual acceptance is set to $\epsilon_{tol}=C \Delta t ^{p+1}$, where $C$ is a user-specified constant and $p$ is the order of convergence for the temporal discretization scheme. {This choice of tolerance ensures that residual error is smaller than the LTE of time integration scheme utilized.}  If the residual norm is not small enough compared with $\epsilon_{tol}$, Krylov subspaces will be further augmented in Step K1 in an iterative fashion, until the prescribed tolerance is reached. 
After the solution is accepted, the resulting matrix $S_1$ is usually dense. We diagonalize it by performing a reduced SVD decomposition of $S_1$ to ensure that only the dominant singular values and essential singular vectors are retained. 

    %
    %

\item [Step K5.] {\em Null-space correction via a LoMaC projection \cite{guo2023local}.} 
In time-dependent problems, the steady-state solution is characterized by the null space of the diffusion operator. 
Accounting for the null space in the solution construction ensures that the solution converges to the correct equilibrium states, leading to a significant improvement in accuracy. 
{To ensure accurate representation of the null-space components, we propose to apply the LoMaC approach \cite{guo2023local},
which is designed to correct the numerical solution so that its projection to the null space of the operator is consistent with that of the reduced macroscopic model.} 
In particular, we perform the decomposition $U_1S_1{V_1}^T={\Tilde{\bF}}_1+\bF_2$, 
where, ${\Tilde{\bF}}_1$ is the projection of the solution onto the null space and $\bF_2$ is the remainder. In the case of heat equation with periodic boundary conditions, the null space is the constant mode vector, ${\Tilde{\bF}}_1= \tilde{n} \mathbf{1} \otimes \mathbf{1}^T$, with $\tilde{n} = \Delta v_1 \Delta v_2 (\mathbf{1}^T U_0) S_0 (V_0^T \mathbf{1})$ being the average density of the solution. 
The component $\bF_2 = U_1S_1V_1^T - \Tilde{\bF}_1=[U_1,{\bf 1}]\text{\bf diag}(S_1,-\Tilde{n})[V_1,{\bf{1}}]^T$ contains zero total mass.

The subsequent step involves adjusting the solution's mass to match the initial condition's mass $n$ and then
apply a truncation to $f_2$, retaining only the most significant singular vectors, denoted by $\mathcal{T}_{\epsilon}(f_2)$ (see Algorithm \ref{alg:trunc}). Note that $\tilde{n}$ does not necessarily match exactly with the total mass $n$, due to discretization errors. We adjust the solution as:
\begin{align}
    n \mathbf{1} \otimes \mathbf{1}^T + \mathcal{T}_{\epsilon}(f_2), 
\end{align}
where we use $n$ in place of $\tilde{n}$ to project the low-rank solution onto the null space with correct total mass, and $\mathcal{T}_{\epsilon}(f_2)$ contains zero total mass and only significant modes to optimize efficiency. 
\end{enumerate}

\begin{figure}
    \centering
    \includegraphics[width=0.9\textwidth]{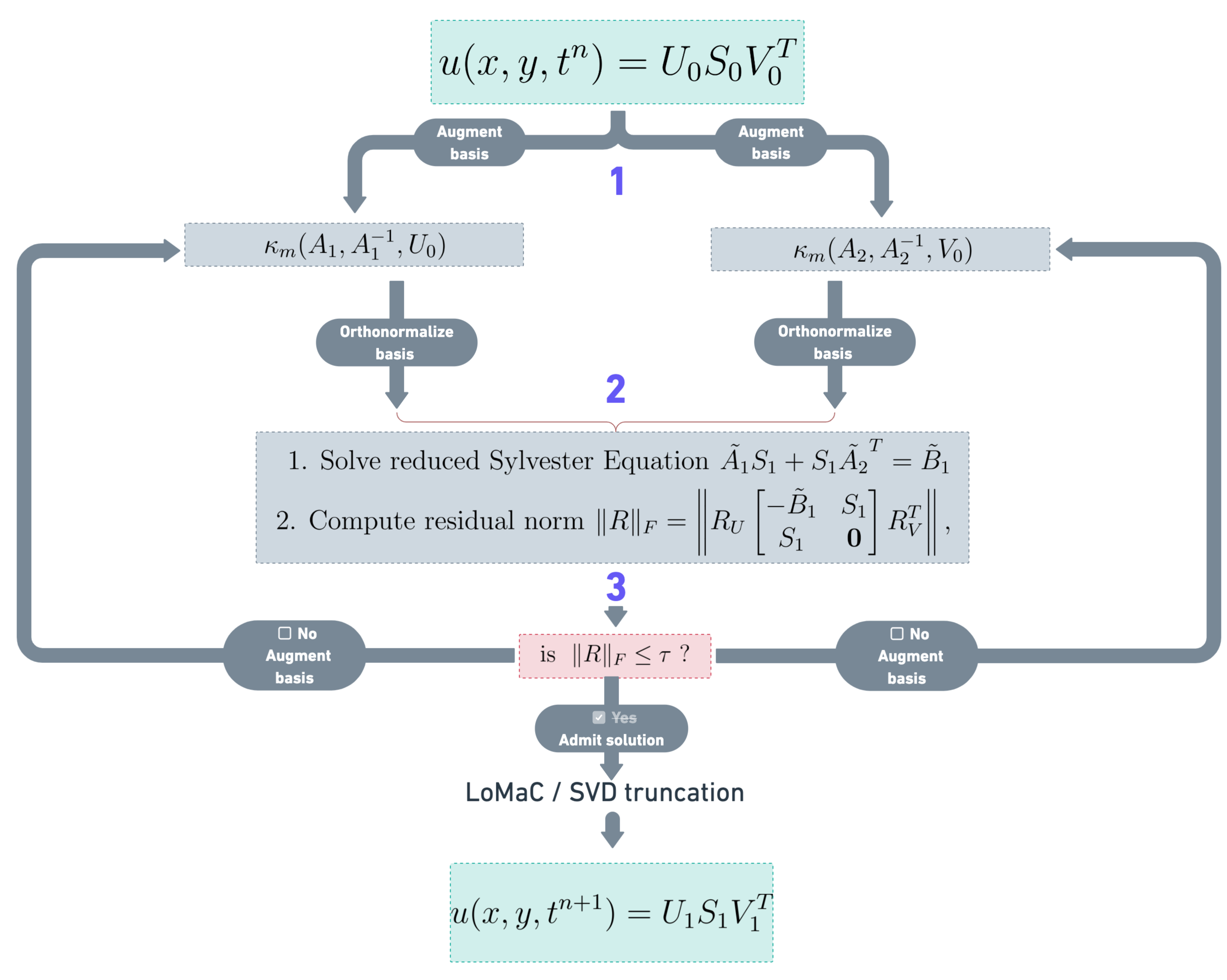}
    \caption{Flow-chart of the extended-Krylov-based implicit adaptive-rank algorithm with the LoMaC projection. }
    \label{flow-chart}
\end{figure}

We summarize the whole procedure in Algorithm \ref{alg:BE-LR}, and illustrate the algorithm flowchart in Figure \ref{flow-chart}.

\begin{algorithm}[t]
\caption{\label{alg:trunc}Truncation Procedure, $\mathcal{T}_\epsilon$}
\SetAlgoNlRelativeSize{-2}
\SetNlSty{textbf}{(}{)}
\KwIn{Bases $U$, $V$; matrix of coefficients $S$; tolerance $\epsilon$.}
\KwOut{Truncated bases $\Tilde{U}$, $\Tilde{V}$; truncated singular values $\Tilde{S}$.}

\If{Bases $U$ and $V$ are not orthonormal}{
    Perform reduced QR decomposition: $Q_1R_1 = U$, $Q_2R_2 = V$\;
    Update left singular vectors: $U \leftarrow Q_1$\;
    Update right singular vectors: $V \leftarrow Q_2$\;
    Update matrix of singular values: $S \leftarrow R_1 S R_2$\;
}
Perform reduced SVD decomposition: $T_1 \Tilde{S} T_2 = \text{SVD}(S)$, where $\Tilde{S} = \text{diag}(\sigma_j)$\label{alg1:line:svd}\;

Identify the last index $\Tilde{r}$ such that $\sigma_{\Tilde{r}} > \epsilon$\;

Update matrix of singular values: $\Tilde{S} \leftarrow \Tilde{S}(1:\Tilde{r}, 1:\Tilde{r})$\;

Update left singular vectors: $\Tilde{U} \leftarrow U T_1(:, 1:\Tilde{r})$\label{alg1:line:mm1}\;

Update right singular vectors: $\Tilde{V} \leftarrow V T_2(:, 1:\Tilde{r})$\label{alg1:line:mm2}\;

\end{algorithm}

\begin{algorithm}[H]
\caption{Backward Euler Adaptive-Rank Integrator}
\label{alg:BE-LR}
\SetAlgoNlRelativeSize{-2}
\SetNlSty{textbf}{(}{)}
\tcp*[h]{This algorithm solves $A_1 \mathbf{F} + \mathbf{F} A_2^T = \mathbf{B}$ for $\mathbf{F}=U_1S_1V_1^T$, where $\mathbf{B}=U_0S_0V_0^T$.}\\
\KwIn{Initial condition $U_0$, $V_0$, $S_0$; Operators $A_1$, $A_2$; Tolerances $\epsilon$, $\epsilon_{\text{tol}}$; Maximum iterations.}
\KwOut{Updated bases $U_1$, $V_1$; Truncated singular values $S_1$.}

\For{m = 1 to maximum iterations}{
   \tcp*[h]{Step K1.}\\
   $U, R \leftarrow \texttt{qr}(\kappa_m(A_1, A_1^{-1}, U_0))$\label{alg:line:qr1}\;
   $V, R \leftarrow \texttt{qr}(\kappa_m(A_2, A_2^{-1}, V_0))$\label{alg:line:qr2}\; 
   Compute $R_U = \texttt{qr}([U, A_1U])$ and $R_V = \texttt{qr}([V, A_2V])$\label{alg:line:qr3}\;
   \tcp*[h]{Step K2.}\\
   Compute $\tilde{A}_1$, $\tilde{A}_2$, and $\tilde{B}_1$ from Equations \eqref{eq: tA1} and \eqref{eq: bB1}\label{alg:line:matrices}\;
   Solve reduced Sylvester equation $\tilde{A}_1 S_1 + S_1 \tilde{A}_2^T = \tilde{B}_1$\label{alg:line:sylv}\; 
   \tcp*[h]{Step K3.}\\
   Compute the residual $\|\mathbf{R}\| = \left\| R_U \begin{bmatrix} -\tilde{B}_1 & S_1 \\ S_1 & 0 \end{bmatrix} R_V^T \right\|$\label{alg:line:residual}\; 

   \If{$\|\mathbf{R}\| \geq \epsilon_{\text{tol}}$}{
      Reject solution and return to Step K1 to augment bases further\;
   }   
   \Else{
      \If{operator does not contain a null-space}{
         \tcp*[h]{Step K4.}\\
         Truncate $\mathcal{T}_\epsilon(U_1 S_1 V_1^T)$\label{alg:line:trunc}\;
         break\;
      }
      \Else{
         \tcp*[h]{Step K5.}\\
         \textbf{Input:} Moments to be conserved, solution bases $U_1$, $V_1$, $S_1$, and truncation threshold $\epsilon$\label{alg:line:lomac1}\;
         \textbf{Output:} Solution with corrected moments, updated $U_1$, $V_1$, $S_1$\label{alg:line:lomac2}\;
         break\;
      } 
      Exit loop\;
   }
}
\end{algorithm}

%
%

\subsection{Computational complexity analysis of the extended-Krylov implicit adaptive-rank algorithm}
\label{sec: 2.3}

In this subsection, we analyze the computational complexity of Algorithm \ref{alg:BE-LR} for backward Euler.  For conciseness, we consider the same spatial resolution per dimension, i.e., \( N=N_1=N_2 \). 
Starting from an initial condition \( U_0S_0V_0^T \), where \( U_0 \in \mathbb{R}^{N \times r} \), \( V_0 \in \mathbb{R}^{N \times r} \), \( S_0 \in \mathbb{R}^{r \times r} \), the \( m \)-th extended Krylov iteration constructs a basis of \( U_1 \in \mathbb{R}^{N \times r_m} \), \( V_1 \in \mathbb{R}^{N \times r_m} \), and \( S_1 \in \mathbb{R}^{r_m\times r_m} \), where $r_m=(2m+1)r$.
The computational complexity for each step is estimated as follows:
%
%
\paragraph{Step K1}
Constructing the extended Krylov basis as outlined in lines \eqref{alg:line:qr1} and \eqref{alg:line:qr2} of Algorithm \ref{alg:BE-LR} requires performing the operation $A_1^{-1}U_0$, prompting the need to solve the system $A_1X=U_0$. This is generaly performed via an  LU factorization of $A_1$.  With 1D finite differences, $A_1$ features a tridiagonal structure, which can be factorized optimally with techniques like the Thomas algorithm, with $\mathcal{O}(N)$ complexity per column of $X$  \cite{quarteroni2006numerical}, resulting in an overall complexity of $\mathcal{O}(Nr)$. Additionally, the application of $A_1 U_0$ can also be accomplished with $\mathcal{O}(Nr)$ computational complexity for sparse tridiagonal matrices. 
In regards to orthonormalization, utilizing a modified Gram-Schmidt method for computing the reduced QR decomposition of \( U_1 \) and \( V_1 \), the complexities of lines \eqref{alg:line:qr1}, \eqref{alg:line:qr2} and \eqref{alg:line:qr3}  in Alg. \ref{alg:BE-LR} scale as \( \mathcal{O}\left(N r_m^2\right) \), where in general we expect \( r_m \ll N \), when the solution displays a low-rank structure. 

%
%

\paragraph{Step K2}
For sparse matrices \( \Tilde{A}_1 \) and \( \Tilde{A}_2 \), line \eqref{alg:line:matrices} involves matrix-matrix multiplications with complexities of \( \mathcal{O}\left(Nr_m^2+N r_m\right) \), while \( \Tilde{B}_1 \) results in complexities of \( \mathcal{O}\left( N r_m r+ r_m r^2+ r_m^2 r \right) \). The solution of the Sylvester equation in line \eqref{alg:line:sylv} results in complexities of \( \mathcal{O}(r_m^3) \) when employing the Bartels–Stewart algorithm \cite{golub2013matrix}.

%
%

\paragraph{Step K3}
The residual computation in line \eqref{alg:line:residual} is of \( \mathcal{O}(r_m^3) \) (i.e., independent of $N$) due to the matrix multiplications, with an additional \( \mathcal{O}(r_m^2) \) for the Frobenius norm.

%
%

\paragraph{Step K4}
The SVD decomposition of the small-sized matrix \( {S_1} \) in line \eqref{alg1:line:svd} of Algorithm \ref{alg:trunc}  is of \( \mathcal{O}(r_m^3) \), also independent of $N$. Updating the bases in lines \eqref{alg1:line:mm1} and \eqref{alg1:line:mm2} of Algorithm \ref{alg:trunc} requires \( \mathcal{O}( Nr_m\Tilde{r}) \).

%
%

\paragraph{Step K5}
The dominant cost of the LoMaC projection revolves around the truncation of $\bF_2$. The QR decomposition 
and matrix multiplications of bases $[U_1,{\bf 1}]$ and $[V_1,{\bf 1}]$ is of \( \mathcal{O}(N(r_m +1 )^2) \).
The matrix multiplication and SVD decomposition of the small-sized matrix \( \text{\bf diag}(S_1,-\Tilde{n}) \)  is of \( \mathcal{O}((r_m +1 )^3) \).
\\
\\
The computational complexity scalings are summarized in Table \ref{tab:complexities}, where we conclude that the overall algorithm should scale as $\mathcal{{O}}(N ({r_m+1})^2)$, which motivates keeping the overall rank $r$ as low as possible.

\begin{table}[t]
\centering
\caption{\label{tab:complexities}Computational complexities of Algorithm \ref{alg:BE-LR}.}
\begin{tabular}{|l|l|}
\hline
\textbf{Line(s)} & \textbf{Complexity}    \\ \hline
\eqref{alg:line:qr1}, \eqref{alg:line:qr2}     & $\mathcal{O}\left( N(r +r_m^2)\right)$      \\ \hline
\eqref{alg:line:qr3}           & $\mathcal{O}\left( Nr_m^2 \right)$      \\ \hline
\eqref{alg:line:matrices}          & $\mathcal{O}\left(Nr_m^2+N r_m\right) $  and $\mathcal{O}\left( N r_m r+ r_m r^2+ r_m^2 r \right)$     \\ \hline 
\eqref{alg:line:sylv}  & $\mathcal{O}(r_m^3) $     \\ \hline
\eqref{alg:line:residual}           & $\mathcal{O}(r_m^3)$     \\ \hline
\eqref{alg:line:residual} (Frobenius norm) & $\mathcal{O}(r_m^2)$ \\ \hline
\eqref{alg:line:trunc}  (Algorithm \ref{alg:trunc})       & $ \mathcal{O}(Nr_m\Tilde{r}+ (r_m)^3)$ \\ \hline
{\eqref{alg:line:lomac1}-\eqref{alg:line:lomac2} -- Step K5}     & $\mathcal{{O}}(N ({r_m+1})^2 + (r_m+1)^3)$ \\ \hline 
\end{tabular}
\end{table}

%
%

\subsection{Krylov-based adaptive-rank high-order DIRK integrators}
\label{sec: 2.4}

In this section, we extend the adaptive-rank backward Euler implicit integrator to  high-order DIRK time discretizations. DIRK methods can be expressed with the following Butcher table presented in Table \ref{table:DIRK}.
\begin{table}[h!]
\centering
\caption{Butcher table for an $s$ stage DIRK scheme. Here $s$ represents the number of stages in DIRK, $c_i$ represents the intermediate stage at which the solution is being approximated, $a_{ij}$ is a lower-triangular matrix with coefficients used to approximate the solutions at intermediate stages, and $b_j$ represents the quadrature weights to update DIRK solution in the final step.}
\label{table:DIRK}
\begin{tabular}{c|llll}
    $c_1$&$a_{11}$&0&$\hdots$&0\\
    $c_2$&$a_{21}$&$a_{22}$&$\hdots$&0\\
    $\vdots$&$\vdots$&$\vdots$&$\ddots$&$\vdots$\\
    $c_s$&$a_{s1}$&$a_{s2}$&$\hdots$&$a_{ss}$\\
    \hline
    &$b_1$&$b_2$&$\hdots$&$b_s$
    \end{tabular}
\end{table}

We consider stiffly accurate DIRK methods, where the coefficients $b_i= a_{si}$ for $i=1,\cdots s$. Thus the solution in the final update is the solution in the final DIRK stage. The corresponding DIRK scheme for the matrix differential equation \eqref{Semi-Discrete} from $t^{(n)}$ to $t^{(n+1)}$ can be written as follows:
\begin{subequations}
\label{eq:RK}
\begin{align}
\begin{split}\label{eq:RKstage_a}
    \mathbf{F}^{(k)} &= \mathbf{F}^n + \Delta t\sum\limits_{\ell=1}^{k}{a_{k\ell}\mathbf{Y}_{\ell}},\qquad k=1,2,...,s.
\end{split}
\\
\begin{split}
    \mathbf{Y}_k &= \mathbf{D}(\mathbf{F}^{(k)};t^{(k)}), \qquad t^{(k)} = t^n + c_k \Delta t, \qquad k=1,2,...,s, 
\end{split}
\\
\begin{split}\label{eq:RK_finalstep}
    {\bf{F}}^{n+1} &= \mathbf{F}^{(s)} = {\bf{F}}^{n} + \Delta t\sum\limits_{k=1}^{s}{b_k\mathbf{Y}_k} .
\end{split}
\end{align}
\end{subequations}
The DIRK method features a stage-by-stage backward-Euler-like implicit solver, with explicit evaluation of RHS terms from previous RK stages \cite{guo2022low}. In particular, for each $k^{th}$ RK stage, we have
    \begin{equation}
      A_1^{(k)} {\bF^{(k)}} + {\bF}^{(k)} A_2^{(k)} = {\bf B}^{(k)}, 
        \label{eq: Syl_k}
    \end{equation}
    with 
    \begin{eqnarray}
 && A_1^{(k)} =\frac{1}{2}I - \Delta t a_{kk} D^{(k)}_1, \qquad A_2^{(k)} = \frac{1}{2}I - \Delta t a_{kk} D^{(k)}_2, \qquad
       {\bf B}^{(k)} = \mathbf{F}^n +  \Delta t\sum\limits_{\ell=1}^{k-1}{a_{k\ell}\mathbf{Y}_{\ell}}
       \label{disc-sylv}
    \end{eqnarray} 
Here, $\bB^{(k)}$ can be explicitly evaluated from solutions at previous RK stages, owing to the Kronecker product structure of the differential operator $\mathbf{D}(\mathbf{F}^{(k)};t^{(k)})$, e.g. Eq. \eqref{eq:FP}. 
We employ our proposed low-rank Krylov-based Sylvester solver in Section~\ref{sec: 2.1} to solve \eqref{eq: Syl_k} stage-by-stage. 
To further improve algorithmic efficiency, with the consideration that the solution spaces across DIRK stages will be highly correlated, we construct the set of Krylov subspaces for the initial DIRK stage and use it throughout all DIRK stages unless the specified error tolerance is not satisfied. Specific steps are as follows: 
\begin{enumerate}
    \item [Step 1.] Prediction of Krylov basis functions: construct a set of orthonormal bases $U_1$ and $V_1$ from the Krylov-based low-rank implicit solver at the first DIRK stage, i.e. backward Euler with time stepping size $c_1 \Delta t$.
    \item [Step 2.] For $k=1:s$ (per DIRK stage)
    \begin{enumerate}
        \item Solve reduced Sylvester equation for $S^{(k)}$,
    \begin{equation}
     \tilde{A}^{(k)}_1 S^{(k)} + S^{(k)} \tilde{A}^{(k)}_2 = \tilde{B}^{(k)}, 
        \label{eq: Sk}
    \end{equation}
    with \[
    \tilde{B}^{(k)} =  U_1^T \mathbf{F}^n V_1 +  \Delta t\sum\limits_{\ell=1}^{k-1}{a_{k\ell}\tilde{Y_{\ell}}}, 
    \]
    Here $\tilde{A}^{(k)}_1 = \frac12 I - a_{kk} \Delta t \tilde{D}_1$ as in \eqref{eq: tA1} and similarly for $\tilde{A}^{(k)}_2$. $\tilde{Y_{\ell}} = U_1^T \mathbf{Y}_{\ell} V_1 \in \mathbb{\bR}^{r \times r}$. Solve for $S^{(k)}$ and obtain the intermediate RK solutions as $U_1 \bS^{(k)} V_1^T$. 
To further improve computational efficiency, from \eqref{eq: Sk}, we have
\begin{align*}
    \tilde{Y}_{\ell} &= \frac{1}{a_{\ell \ell}\Delta t}\left( {S^{(\ell)}} -  {\tilde{B}^{(\ell)}}\right),
\end{align*}
leading to an efficient computation of $\tilde{B}^{(k)}$ 
\begin{equation}
     \tilde{B}^{(k)} =  \tilde{B}_1 + \Delta t \sum_{\ell=1}^{k-1} \frac{a_{k\ell}}{a_{\ell \ell}} \left( {S^{(\ell)}} -  {\tilde{B}^{(\ell)}}\right),
\end{equation}
    with $\tilde{B}_1$ defined in \eqref{eq: bB1}. This avoids the need of evaluating the full size $\bB^{(k)}$ in \eqref{eq: Syl_k}.
    \item Efficiently evaluate the residual norm, 
    \[
 \left\|{\bR}^{(k)} \right\|= \left\|  {R}_U
\begin{bmatrix}
  - \tilde{B}^{(k)} &  {S}^{(k)} \\
   {S}^{(k)} & 0
\end{bmatrix}
 {R}_V^T \right\|,
\]
where $R_U$ and $R_V$ are upper triangular matrices from a reduced QR decomposition of $[U_1, A_1 U_1]$ and $[V_1, A_2 V_1]$ respectively, which can be done once for all DIRK stages. 
\item Adaptive criteria via a prescribed tolerance: we compare the computed residual norm $\|R\|$ with the given error tolerance $\epsilon^{(k)}_{tol}$. If the tolerance is met, then move to the next DIRK stage with $\bF^{(k)} = U_1 S^{(k)} V_1^T$; otherwise we go back to Step 1 to further augment Krylov subspaces. 
\end{enumerate} 
\item [Step 3.] If the operator contains a null-space, perform LoMaC projection on constructed solution (See step K5. of section \ref{sec: 2.2}); otherwise, 
 perform truncation on evolved solution, i.e. $\mathcal{T}_\epsilon(U_1 S_1 V_1^T)$.
    
\end{enumerate}
The proposed adaptive-rank DIRK integrator is summarized in Algorithm \ref{DIRK-LR}. 
\begin{algorithm}[!ht]
\caption{s-Stages Adaptive-Rank DIRK Integrator}
\label{DIRK-LR}
\SetAlgoNlRelativeSize{-2}
\SetNlSty{textbf}{(}{)}
\tcp*[h]{This algorithm evolves the adaptive-rank solution by a high order DIRK integrator. }\\
\KwIn{Initial condition $U_0$, $V_0$, $S_0$; Operators $\{A^{(1)}_1, \cdots A^{(s)}_1\}$, $ \{A^{(1)}_2, \cdots, A^{(s)}_2\}$; Butcher table $\{a_{ij}\}$; Time step size $\Delta t$; Tolerances $\epsilon$, $\vec{\epsilon}_{\text{tol}}$; Maximum iterations.}
\KwOut{Updated bases $U_1$, $V_1$; Truncated singular values $S_1$.}
\For{m = 1 to maximum iterations}{
  \tcp*[h]{Step 1.}\\
  $U_1, R \leftarrow \texttt{qr}(\kappa_m(A^{(1)}_1, A_1^{-1}, U_0))$\;
  $V_1, R \leftarrow \texttt{qr}(\kappa_m(A^{(1)}_2, A_2^{-1}, V_0))$\;
  Compute $\tilde{B}^{(1)}$\;
  \tcp*[h]{Step 2.}\\
  \For{$k = 1$ to $s$}{
    Compute $\tilde{B}^{(k)} = \tilde{B}^{(1)} + \Delta t \sum_{\ell=1}^{k-1} \frac{a_{k\ell}}{a_{\ell \ell}}\left( S^{(\ell)} - \tilde{B}^{(\ell)} \right)$\;
    Solve Sylvester equation $\tilde{A}^{(k)}_1 S^{(k)} + S^{(k)} \tilde{A}^{(k)}_2 = \tilde{B}^{(k)}$\;
    Compute $R_U = \texttt{qr}([U_1, A_1^{(k)} U_1])$ and $R_V = \texttt{qr}([V_1, A_2^{(k)} V_1])$\;
    Compute the residual $\|\mathbf{R}^{(k)}\| = \left\| R_U \begin{bmatrix} -\tilde{B}^{(k)} & S^{(k)} \\ S^{(k)} & 0 \end{bmatrix} R_V^T \right\|$\;
    \If{$\|\mathbf{R}^{(k)}\| \geq \epsilon^{(k)}_{\text{tol}}$}{
      Reject step and return to Step 1 to augment bases further\;
    }
    \Else{
      Compute and store $\frac{1}{a_{kk}}\left( S^{(k)} - \tilde{B}^{(k)} \right)$ and proceed to the next stage\;
    }
  }
  \tcp*[h]{Step 3.}\\
  \If{all stages are accepted}{
    \If{operator does not contain a null-space}{
      Truncate $\mathcal{T}_\epsilon(U_1 S_1 V_1^T)$\;
      break\;
    }
    \Else{
      \tcp*[h]{Perform LoMaC correction.}\\
       LoMaC correction \\
     \textbf{Input:} Moments to be conserved, solution bases $U_1$, $V_1$, $S_1$, and truncation threshold $\epsilon$\;
      \textbf{Output:} Solution with corrected moments, updated $U_1$, $V_1$, $S_1$\;
      break\;
    }
  }
}
\end{algorithm}

%
%
\section{Adaptive-rank implicit integrators for the multi-species nonlinear Fokker-Planck equation}
\label{sec:fp}

\subsection{Multi-species LBFP model}
The multi-species LBFP equation is a simplified collisional model describing the collisional relaxation of multiple plasma species. In this study, we employ a recently proposed formulation of the LBFP model \cite{filbet2022fokker} that features exact conservation properties and the H-theorem for entropy dissipation, which for species $\alpha$ reads:
\begin{equation}
    \label{eq:lbfp}
    \frac{\partial f_{\alpha}}{\partial t} = 
    \sum^{N_s}_{{\beta=1}} \nu_{\alpha\beta} \nabla_v \cdot \left[
    D_{\alpha\beta} \nabla_v f_{\alpha} + \left(\vec{v} - \vec{u}_{\alpha\beta} \right) f_{\alpha}
    \right].
\end{equation}
Here, most quantities are as defined earlier, but without loss of generality, we consider a two-dimensional velocity space $\vec{v} = \left\{v_1,v_2\right\} \in {\mathbb R}^2$ for simplicity. Following \cite{filbet2022fokker}, the coefficients are defined as follows: $D_{\alpha\beta}$ is the diffusion coefficient (related to the temperature, see below), $\vec{u}_{\alpha} = \frac{\vec{\gamma}_{\alpha}}{n_{\alpha}}$ is the drift velocity,  $\vec{u}_{\alpha\beta} = \frac{\vec{u}_{\alpha}+\vec{u}_{\beta}}{2}$ is the mixed drift velocity. 
Here $n_{\alpha} = \left< 1, f_{\alpha} \right>_v$ is the number density, $\vec{\gamma}_{\alpha} = \left< \vec{v}, f_{\alpha} \right>_v$ is the particle flux, and $\left< A, B \right>_v = \int_{{\mathbb R}^2} d^2 v AB$ is a shorthand notation for the velocity space inner-product between functions $A$ and $B$. We consider proton and electron species, $\alpha \in \left\{ p, e\right\}$, with collisional coefficients given by:
\begin{eqnarray*}
\nu_{\alpha\beta}={2^{5/2}}e_\alpha^2 e_\beta^2 n_\beta \frac{ m_\beta}{m_\alpha+ m_\beta}\frac{1}{(v_{th_{\alpha}} +v_{th_{\beta}})^{\frac{3}{2}}}, \; v_{th_{\alpha}}=\sqrt{\frac{T_\alpha}{m_\alpha}}, \label{param1}\\
D_{\alpha\beta}=\frac{T_{\alpha,\beta}}{m_\alpha}, \; T_{\alpha\beta}= \frac{m_\alpha T_\beta-m_\beta T_\alpha}{m_\alpha+m_\beta}+\frac{m_\alpha m_\beta}{m_\alpha + m_\beta}|\vec{u}_\beta -\vec
{u}_\alpha|^2. \label{param2}
\end{eqnarray*}

The macroscopic conservation laws for mass, momentum, and energy are obtained for species $\alpha$ by projecting Eq. \eqref{eq:lbfp} onto the $\vec{\phi} = \left\{ 1, \vec{v}, \frac{1}{2}Tr\left(\vec{v}\vec{v}\right) \right\}$ subspace, to find:
\begin{eqnarray}
    \partial_t n_{\alpha} = 0, \label{continuityODE}\\
    \partial_t \vec{\gamma}_{\alpha} = \frac{1}{2}\sum_{\beta\neq\alpha}^{N_s}\nu_{\alpha\beta} n_{\alpha} \left(\vec{u}_{\alpha} - \vec{u}_{\beta} \right), \label{momentumODE}\\
    \partial_t {\cal E}_{\alpha} = \sum^{N_s}_{\beta\neq\alpha} \nu_{\alpha\beta} \left( 2D_{\alpha\beta} n_{\alpha} - 2{\cal E}_{\alpha} + \frac{1}{2} \vec{\gamma}_{\alpha} \cdot (\vec{u}_{\alpha}+\vec{u}_{\beta}) \right).\label{energyODE}
\end{eqnarray}
Here, ${\cal E}_{\alpha} = \frac{1}{2}\left< Tr\left(\vec{v}\vec{v}\right), f_{\alpha}\right>_v = \frac{\vec{u}_{\alpha}\cdot\vec{\gamma}_{\alpha}}{2} + \frac{d}{2}\frac{n_{\alpha}T_{\alpha}}{m_{\alpha}}$ is the specific total energy density. The macroscopic conservation theorems can readily be shown by multiplying each species conservation equation with their respective mass and summing them over to yield \cite{filbet2022fokker}:
\begin{eqnarray}
    \sum^{N_s}_{\alpha} m_{\alpha} \partial_t n_{\alpha} = 0, \qquad
    \sum^{N_s}_{\alpha} m_{\alpha} \partial_t \vec{\gamma}_{\alpha} = \vec{0}, \qquad
    \sum^{N_s}_{\alpha} m_{\alpha} \partial_t {\cal E}_{\alpha} = 0.
\end{eqnarray}

The LBFP system also satisfies the Boltzmann ${\cal H}$-theorem \cite{filbet2022fokker},
\begin{equation}
    \label{eq:h_theorem}
    \frac{d{\cal H}}{dt} \le 0.
\end{equation}
Here, ${\cal H}\left[f\right] = \sum^{N_s}_{\alpha} \left< f_{\alpha},\ln f_{\alpha} \right>_v$ is the total entropy functional and Eq. \eqref{eq:h_theorem} can be shown to monotonically decay until $f_{\alpha} = \bar{f}^M_{\alpha}\left( \vec{v}; \vec{\bar{u}}, \bar{T} \right) = \frac{n_{\alpha}}{2\pi\bar{T}/m_{\alpha}}\exp\left( -\frac{m_{\alpha}}{2\bar{T}} \left| \vec{v} - \vec{\bar{u}} \right|^2 \right)$, with $\vec{\bar{u}}$ and $\bar{T}$ the equilibrium drift and temperature of the system defined from the momentum and energy conservation theorems:
\begin{eqnarray*}
    \vec{\bar{u}} = \frac{\sum^{N_s}_{\alpha} m_{\alpha} n_{\alpha} \vec{u}_{\alpha}\left(t\right)}{\sum^{N_s}_{\alpha} m_{\alpha} n_{\alpha}} \\
    \bar{T} = \frac{\frac{1}{2} \sum^{N_s}_{\alpha} m_{\alpha} n_{\alpha} u^2_{\alpha} \left(t\right) + \sum^{N_s}_{\alpha} n_{\alpha} T_{\alpha}\left(t \right) - \frac{\bar{u}^2}{2} \sum^{N_s}_{\alpha} m_{\alpha} n_{\alpha}}{\sum^{N_s}_{\alpha} n_{\alpha}}.\label{eq:Temp-eq}
\end{eqnarray*}
\subsection{Temporal update of the LBFP model}
\label{sec:3.2}
To discretize the LBFP equation for a given species {$\alpha \in  \{i,e\}$}, we consider a two-dimensional tensor grid. This grid is characterized by $N_{v}$ uniform discretization points in each dimension, supporting the initial distribution ${\bf {\mathbf F}}_{\alpha}$. Here, ${\bf {\mathbf F}}_{\alpha}$ denotes the discrete solution matrix that encapsulates the initial conditions of the distributions $f_\alpha$ on our tensor grid. The velocity domain is dimensioned 
based on the thermal velocity of each species, defined as $v_{th,\alpha}=\sqrt{\frac{T_\alpha}{m_\alpha}}$, as follows: for each species $\alpha$, we set the spatial domain to span from $-10v_{th,\alpha}$ to $10v_{th,\alpha}$ in both the $v_1$ and $v_2$ directions. This methodical discretization ensures that our distributions are well-supported within the specified domain. 

Next, a truncated SVD is applied to the discrete initial condition matrix ${\bf {\mathbf F}}_{\alpha}$. The result of this decomposition is expressed as ${\bf {\mathbf F}}_{\alpha} = U_{\alpha,0} S_{\alpha,0} V_{\alpha,0}^T$, where $U_{\alpha,0} \in \mathbb{R}^{N_{v} \times r_\alpha}$, $V_{\alpha,0} \in \mathbb{R}^{N_{v} \times r_\alpha}$, and $S_{\alpha,0} \in \mathbb{R}^{r_\alpha \times r_\alpha}$. It is important to acknowledge that, unlike the heat equation, the LBFP equation is nonlinear due to the dependence of the coefficients on the moments of the solution. We linearize it by evolving the moment equations \eqref{continuityODE}-\eqref{energyODE} in time (as proposed in \cite{taitano2015charge}, but restricted here to the collisional terms) with {the same temporal integrator as the adaptive-rank integrator for the kinetic equation.} 
Subsequently, we construct our differentiation matrices for the resulting linearized problem and proceed with our low-rank integration. As noted earlier, the evolution of the distribution function might not maintain the moments due to numerical discretizations and SVD truncation. To mitigate this issue, we employ a LoMaC post-processing technique \cite{guo2022conservative} (also see section \ref{sec:3.3}). This approach adeptly corrects the moments of our distribution while simultaneously preserving the rank-adaptive nature of our algorithm.

To evolve the solution from time $t^{(n)}$ to $t^{(n+1)}$, we use an s-stage DIRK scheme as outlined in Section \ref{sec: 2.4}. We follow the steps below: 
\begin{enumerate}
    \item \textbf{Step 1: {Integrating the macroscopic ODEs.}}
    First, we integrate the ordinary differential equations (ODEs) given by equations \eqref{continuityODE}-\eqref{energyODE} implicitly in time using the same DIRK scheme used to evolve the distributions. This results in a set of nonlinear equations related to the moments of species $\alpha$. We solve these equations using Newton's method from {$t^n$ to $t^{n+1}$}. This step allows us to compute the values of moments $n_{\alpha}^{(k)}$, $\Vec{\gamma}_{\alpha}^{(k)}$, and ${\cal E}_{\alpha}^{(k)}$ (which represent the mass, momentum, and energy for each species) at different DIRK stages (indexed by $k=1, \cdots, s$). 
    \item \textbf{Step 2: Determining Coefficients and Matrices.}
    Utilizing the moments from Step 1, we then calculate the coefficients in equations \eqref{param1}-\eqref{param2} for the right-hand side (RHS) of the LBFP equation in Eq. \ref{eq:lbfp}, with which we construct the differentiation matrices ${A^{(k)}_{1,\alpha}, A^{(k)}_{2,\alpha}}$. These are similar to the matrices mentioned in equation \eqref{disc-sylv} but have an extra subscript $\alpha$ to indicate they are specific to each species. The construction of these matrices depends on the finite difference approximations of derivatives imposed in our domain. Here, we use an equilibrium-preserving second-order discretization proposed by Chang and Cooper \cite{chang1970practical}. For more information on how these matrices are constructed, see Appendix~\ref{app:chang-cooper}.
    \item \textbf{Step 3: Evolving the Solution.}
    Using the differentiation matrices formed in Step 2, we apply the low-rank integrator algorithm (see Algorithm \ref{DIRK-LR}) to evolve the solution from $t^{(n)}$ to $t^{(n+1)}$. This process evolves the bases $U_{\alpha, 1}$, $V_{\alpha, 1}$, and the matrix $S_{\alpha, 1}$ from $t^{(n)}$ to $t^{(n+1)}$.
    \item \textbf{Step 4: Correcting the Solution.}
    The updated solution ${\mathbf F}^{\ast}_{\alpha}=U_{\alpha, 1}S_{\alpha, 1}V_{\alpha, 1}^T$ may not accurately preserve mass, momentum, and energy (hence the $\ast$ symbol), due to the numerical discretization errors and SVD truncation in previous steps. To fix this, we apply a LoMaC projection \cite{guo2023local} as a post-processing step to correct the macroscopic quantities of the computed solution. We elaborate on this procedure further in the next section.
\end{enumerate}

\subsection{Local Macroscopic Conservation (LoMaC) formulation for LBFP}
\label{sec:3.3}
The LoMaC method \cite{guo2023local} corrects the moments of the computed particle distributions every time step. 
For brevity, we omit the species subscript $\alpha$, with the understanding that the formulation applies to each species $\alpha$.

The LoMaC method corrects for the loss of moment conservation, notably mass, momentum, and energy, during the integration and subsequent SVD truncation steps in the integration of the LBFP equation. As for the LoMaC procedure for the heat equation in Section~\ref{sec: 2.2}, we propose to decompose $\bF = \tilde{\bF}_1 + \bF_2$ and perform moment correction in the following steps.

\begin{itemize}
    \item $\tilde{\bF}_1$ is the projection of $\bF$ onto a lower-dimensional subspace,
    \begin{align}
    \mathcal{L} = \text{span}\left\{ {\bf 1}_{{ v_{ 1}}\otimes{v_2}}, {\bf v_{ 1}} \otimes {\bf 1}_{v_2}, {\bf 1}_{v_{\bf 1}}\otimes {\bf v_2}, {\bf v_{ 1}}^2 \otimes {\bf 1}_{v_2}+ {\bf 1}_{v_{ 1}}\otimes {\bf v_2}^2 \right\},
    \label{subspace}
\end{align}    
for preservation of mass, momentum, and energy densities. To ensure proper decay of $\tilde{\bF}_1$ and respect the equilibrium Maxwellian distribution in velocity directions, we introduce the weighted inner product space with the weight function being the Maxwellian distribution defined by the thermal velocity $v_{\text{th}} = \sqrt{T/m}$, 
\begin{equation}
    {\bf w}^{(1)} = \exp\left(-\frac{v_1^2}{2 v_{\text{th}}^2}\right), \quad {\bf w}^{(2)} = \exp\left(-\frac{v_2^2}{2 v_{\text{th}}^2}\right),
    \label{eq:weights}
\end{equation}
This particular choice
accounts for disparate thermal velocities for different species, and is crucial in preserving physicality of the numerical solution.
With such weight functions, as in \cite{guo2022conservative}, we define the weighted inner product space via 
\begin{equation}
\label{eq: inner_p}
    \langle f, g \rangle_{\mathbf{w^{(1)}}\otimes\mathbf{w^{(2)}}} = \sum_i \sum_j f_{ij} g_{ij} w^{(1)}_{i}w^{(2)}_{j}, \quad \|f\|_{\mathbf{w}} = \sqrt{\langle f, f \rangle_{\mathbf{w}}},
\end{equation}
with $\mathbf{w} \in \mathbb{R}^{N_v}$ and each $w_i = w(v_i) \Delta_v$ signifying the quadrature weights for $v$-integration, employing the weight function $w(v)$. Let $\tilde{\bF}_1 :=  P_{\mathcal{L}}({\mathbf F})$ be an orthogonal projection onto $\mathcal{L}$, with its explicit construction in equation \eqref{eq: PF1} elaborated in Proposition~\ref{prop:1} below. 
    \item $\bF_2 = {\mathbf F}^{\ast} - \tilde{\bF}_1$ lives in the orthogonal complement of subspace $\mathcal{L}$.  We perform truncation on the solution $\bF_2$ using Algorithm \ref{alg:trunc} to obtain $\mathcal{T}_\epsilon\left(\bF_2\right)=\mathcal{U}\mathcal{S}\mathcal{V}^T$.
    \item Moment Correction. The moments of low-rank solutions in Proposition~\ref{prop:1} are subject to numerical discretization and low-rank truncation errors. Here, we perform a correction step to match them with the system's exact moments. This crucial correction step ensures the simulation's fidelity to the moments computed from macroscopic ODE system \eqref{continuityODE}-\eqref{energyODE}, effectively rectifying any discrepancies introduced in earlier stages. In particular, we let 
    \begin{equation}
    {\mathbf F}= {\mathbf F}_{1} + \mathcal{T}_\epsilon({\mathbf F}_2)
    = [\mathbf{\Tilde{V}}_1,\mathcal{U}] \cdot \textbf{diag}\left( S^1, \mathcal{S}\right) \cdot [\mathbf{\Tilde{V}}_2,\mathcal{V}]^T. 
 \end{equation}
Here  ${\mathbf F}_{1} = \mathbf{\Tilde{V}}_1 S^1 \mathbf{\Tilde{V}}_2^T$,
where $S^1 =  \textbf{diag}\left( n, u^1, u^2, \left(\frac{\mathcal{E}}{m} - c_{3,y} \cdot n\right), \left(\frac{\mathcal{E}}{m} - c_{3,x} \cdot n\right) \right)$ with ${n}$, ${u}^1$, ${u}^2$, ${\mathcal{E}}$ obtained from integrating the ODE system Eq. \eqref{continuityODE}-\eqref{energyODE}, and $\mathbf{\Tilde{V}}_1$ and $\mathbf{\Tilde{V}}_2$ as specified in Proposition~\ref{prop:1}.
\end{itemize}

\begin{prop}
\label{prop:1}
(Construction of $\tilde{\bF}_1$)
\begin{equation}
   \tilde{\bF}_1 := P_{\mathcal{L}}({\mathbf F}) 
   = \mathbf{\Tilde{V}}_1 \left[ \textbf{diag}\left( \tilde{n}, \tilde{u}^1, \tilde{u}^2, \left(\frac{\tilde{\mathcal{E}}}{m} - c_{3,y} \cdot \tilde{n}\right), \left(\frac{\tilde{\mathcal{E}}}{m} - c_{3,x} \cdot \tilde{n}\right) \right) \right]_{\Tilde{S}^1} \mathbf{\Tilde{V}}_2^T, 
\label{eq: PF1}
\end{equation}
where \begin{itemize}
    \item $\mathbf{\Tilde{V}}_1$ and $\mathbf{\Tilde{V}}_2$ are orthonormal basis for weighted inner product space $\mathcal{L}$ with explicit construction as follows: 
    \begin{align*}
\mathbf{\Tilde{V}}_1 &= \left[c_{1,1}\cdot \mathbf{w_1} \cdot \mathbf{1}, c_{2,1} \cdot \mathbf{w_1} \cdot \mathbf{v_1}, c_{2,1}\cdot \mathbf{w_1} \cdot \mathbf{1}, c_{1,1}\cdot \mathbf{w_1} \cdot \mathbf{1}, c_{4,1} \cdot \mathbf{w_1} \cdot (\mathbf{v_1}^2 - c_{3,1} \cdot \mathbf{1})\right], \\
\mathbf{\Tilde{V}}_2 &= \left[c_{1,2}\cdot \mathbf{w_2} \cdot \mathbf{1}, c_{1,2}\cdot \mathbf{w_2} \cdot \mathbf{1}, c_{2,2} \cdot \mathbf{w_2} \cdot \mathbf{v_2}, c_{4,2} \cdot \mathbf{w_2} \cdot (\mathbf{v_2}^2 - c_{3,2} \cdot \mathbf{1}), c_{1,2}\cdot \mathbf{w_2} \cdot \mathbf{1}\right], 
\end{align*}
where $c_{1,j}=\|\mathbf{1}_{{v_j}}\|_{\mathbf{w_j}}^2$, $c_{2,j}=\|\mathbf{v_j}\|_{\mathbf{w_j}}^2$, $c_{3,j}=\frac{\langle \mathbf{1}_{v_j}, \mathbf{v_j}^2 \rangle_{\mathbf{w}}}{\langle \mathbf{1}_{v_j},\mathbf{1}_{v_j}\rangle_{\mathbf{w_j}}}$, $c_{4,j}= \|\mathbf{v_j}^2-c_{3,j} \mathbf{1}_{v_j}\|_{\mathbf{w_j}}^2$, for $j\in\{1,2\}$
\item $\Tilde{n}$, $\Tilde{u}^1$, $\Tilde{u}^2$, $\mathcal{\Tilde{E}}$ are macroscopic densities of numerical solutions, which could be computed by an efficient low rank integration:
\begin{align*}
\Tilde{n} &= \Delta { v_1} \Delta { v_2} \mathbf{1}^T {\mathbf F}^{\ast} \mathbf{1} = \Delta { v_1} \Delta { v_2} (\mathbf{1}^T U) S (V^T \mathbf{1}),\\
\Tilde{u}^1 &= \Delta { v_1} \Delta { v_2} \mathbf{v}_1^T {\mathbf F}^{\ast} \mathbf{1} = \Delta { v_1} \Delta { v_2} (\mathbf{v}_1^T U) S (V^T \mathbf{1}),\\ 
\Tilde{u}^2 &= \Delta { v_1} \Delta { v_2} \mathbf{1}^T {\mathbf F}^{\ast} \mathbf{v}_2 = \Delta { v_1} \Delta { v_2} (\mathbf{1}^T U) S (V^T \mathbf{v}_2),\\ 
2 \mathcal{\Tilde{E}} &= \Delta { v_1} \Delta { v_2} \left( \mathbf{v}_1^{2T} {\mathbf F}^{\ast} \mathbf{1}+  \mathbf{1}^T {\mathbf F}^{\ast} \mathbf{v}_2^2 \right) = \Delta { v_1} \Delta { v_2} \left( (\mathbf{v}_1^{2T} U) S (V^T \mathbf{1}) + (\mathbf{1}^T U) S (V^T \mathbf{v}_2^2) \right).
\end{align*}
\end{itemize}
\end{prop}
\begin{proof}
It is straightforward to check that $c_{i,j}$, $i=1,.., 4$, $j=1, 2$ are normalization coefficients ensuring that the constructed basis is orthonormal and that the orthogonal projection in the inner product space defined by \eqref{eq: inner_p} gives \eqref{eq: PF1}.
\end{proof}
We now summarize the proposed implicit low rank LBFP integrator with LoMaC projection in Algorithm \ref{DIRK-LR-FPLB}. 
\begin{algorithm}[!ht]
\caption{LBFP Integrator}
\label{DIRK-LR-FPLB}
\SetAlgoNlRelativeSize{-2}
\SetNlSty{textbf}{(}{)}
\tcp*[h]{This algorithm evolves the solution of species $\alpha$ from $t^{(n)}$ to $t^{(n+1)}$.}\\
\KwIn{Initial conditions $U_{\alpha,0}$, $V_{\alpha,0}$, $S_{\alpha,0}$. Initial parameters: densities $n_{\alpha,0}$, drift velocities $\vec{u}_{\alpha,0}$, temperatures $T_{\alpha,0}$.}
\KwOut{Updated bases $U_{\alpha,1}$, $V_{\alpha,1}$, truncated singular values $S_{\alpha,1}$.}

\While{stepping}{
    \For{each species $\alpha$}{
        \For{$k = 1$ to $s$}{
            Compute and store macroscopic quantities $n_{\alpha}^{(k)}$, $\vec{\gamma}_{\alpha}^{(k)}$, $\mathcal{E}_{\alpha}^{(k})$ from Equations \eqref{continuityODE}-\eqref{energyODE}\;
            Construct and store sparse matrices $A_{\alpha,1}^{(k)}$, $A_{\alpha,2}^{(k)}$ (See Appendix \ref{app:chang-cooper})\;
        }
        Evolve solution using $U_{\alpha}^{(n)}$, $S_{\alpha}^{(n)}$, $V_{\alpha}^{(n)}$ via Algorithm \ref{DIRK-LR}\;
        Perform a post-processing LoMaC update to correct the moments, see Section \ref{sec:3.3}\;
    }
}
\end{algorithm}

%
%
\section{Numerical experiments}

\label{sec:num-exp}
This section presents a series of numerical experiments to demonstrate the efficacy of the proposed Krylov-based implicit adaptive-rank algorithm. All simulations were conducted utilizing MATLAB running on a Macbook Pro equipped with a 2.3 GHz Quad-Core Intel Core i7 processor.
 To benchmark our algorithm's performance, we compare our simulation results against those obtained from a full-rank integrator. The full-rank integrator solutions are computed using MATLAB's internal function $\texttt{Sylvester}$ on the full-mesh to solve the Sylvester equation \eqref{eq: Syl}. 

\subsection{Heat equation}

We consider first the following prototype heat equation:
\begin{equation}
\frac{\partial u}{\partial t}=d_1\frac{\partial^2 u}{\partial x^2}+d_2\frac{\partial^2 u}{\partial y^2}, \qquad 0\leq x\leq 1, \qquad 0\leq y\leq 1,
\label{eq:heat}
\end{equation}
where $d_1=d_2=1/2$ represent the diffusion coefficients. 
We consider a rank-2 initial condition:
\begin{equation*}
u(x, y, 0) = 0.5 \exp\left(-400 \left((x-0.3)^2 + (y-0.35)^2\right)\right) + 0.8 \exp\left(-400 \left((x-0.65)^2 + (y-0.5)^2\right)\right),
\end{equation*}
with periodic boundary conditions. 
The spatial grid is discretized with a number of nodes $N_x=N_y=400$, and the time-step is $\Delta t= \lambda \Delta x^2$, with $\lambda$ ranging from $100$ to $900$. The SVD truncation threshold is set to $10^{-10} \sigma_1$, where $\sigma_1$ is the largest singular value in the simulation. 
The spatial differentiation matrices are constructed 
using a finite difference approximation of the second-order derivatives with periodic boundary conditions, which yields circulant matrices. The matrices \({D^{(k)}_1}\) and \({D^{(k)}_2}\), where $k$ is the DIRK stage, analogous to \eqref{disc-sylv} are given by
\[
{D^{(k)}_1}={D^{(k)}_2} = \begin{bmatrix}
-1 & \frac{1}{2} & 0 & \cdots & 0 & \frac{1}{2} \\
\frac{1}{2} & -1 & \frac{1}{2} & \cdots & 0 & 0 \\
0 & \frac{1}{2} & -1 & \cdots & 0 & 0 \\
\vdots & \vdots & \vdots & \ddots & \vdots & \vdots \\
0 & 0 & 0 & \cdots & -1 & \frac{1}{2} \\
\frac{1}{2} & 0 & 0 & \cdots & \frac{1}{2} & -1
\end{bmatrix}.
\]

\begin{figure}
    \centering
    \includegraphics[width=\textwidth]{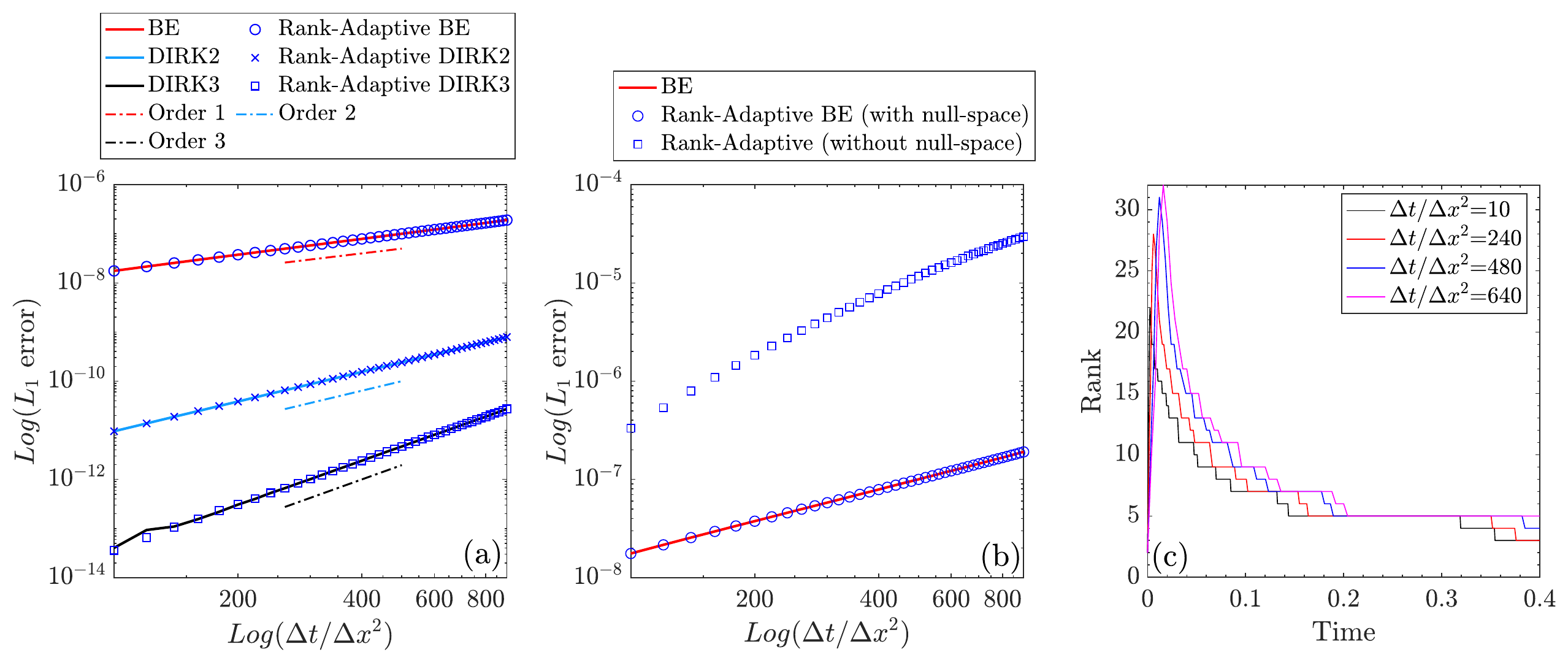}
    \caption{Simulation of the heat equation under periodic boundary conditions employing the BE, DIRK2, and DIRK3 integrators. In Fig. (a), we illustrate a temporal error convergence study spanning a range of $\Delta t/\Delta x^2$ values [100:900], presenting results for both the adaptive-rank integrator and the classical full-rank integrator counterpart. Fig. (b) presents results for both the BE adaptive-rank integrator with and without LoMaC projection, and the classical full-rank integrator counterpart. Fig. (c) shows the rank evolution as a function of time for the first-order adaptive-rank BE integrator with LoMaC projection.
    }
    \label{fig:Null-space-BE}
\end{figure}

In Fig.~\ref{fig:Null-space-BE} (a), 
we assess the $L_1$ error norm by comparing the performance of the adaptive-rank integrator along the LoMaC post-processing scheme against its full-rank counterpart for BE, DIRK2, and DIRK3. The Butcher tables for DIRK2 and DIRK3 are shown in Tables \ref{Butcher:DIRK2} and \ref{Butcher:DIRK3}, respectively. 
The temporal error convergence for our proposed adaptive-rank integrator, depicted by markers, closely matches the performance of the full-rank integrator, represented by solid lines. This comparison reveals that our proposed adaptive-rank algorithm yields solutions with temporal errors comparable to those of the full-rank integrator, while also maintaining a comparatively small solution rank relative to the dimensions of the mesh. 
\begin{table}[h]
\centering
\label{table:ComparisonDIRK}
\begin{minipage}{.4\linewidth}
\centering
\begin{tabular}{c|cc}
 $\gamma$ & $\gamma$ & 0 \\
$1$ &  $1-\gamma$ & $\gamma$ \\
\hline
    & $1-\gamma$    & $\gamma$
\end{tabular} 
\caption{\centering DIRK2 Butcher table with $\gamma=1-\frac{\sqrt{2}}{2}$.}
\label{Butcher:DIRK2}
\end{minipage}%
\begin{minipage}{.7\linewidth}
\centering
\begin{tabular}{c|ccc}
$x$ & $x$ & 0 & 0 \\
$\frac{1+x}{2}$ & $\frac{1-x}{2}$ & $x$ & 0 \\
$1$ & $-\frac{3x^2}{2}+4x-\frac{1}{4}$ & $-\frac{3x^2}{2}-5x+\frac{5}{4}$ & $x$ \\
\hline
 & $-\frac{3x^2}{2}+4x-\frac{1}{4}$ & $-\frac{3x^2}{2}-5x+\frac{5}{4}$ & $x$ 
\end{tabular}
\caption{\centering DIRK3 Butcher table with $x=0.4358665215$.}
\label{Butcher:DIRK3}
\end{minipage}
\end{table}

Next, we utilize the BE adaptive-rank integrator along the LoMaC post-processing scheme for density (zeroth moment) in comparison to the BE adaptive-rank integrator without the LoMaC projection step. Fig. \ref{fig:Null-space-BE} (b) demonstrates that the utilization of LoMaC leads to the correct steady-state solution. In this case, the error convergence of the adaptive-rank integrator, indicated by circle markers, aligns with that of the full-rank integrator, represented by a solid red line. Conversely, when the LoMaC is not utilized, an error is induced due to the non-conservation of the mean of the solution as time progresses. 
This emphasizes the critical role played by LoMaC projection to ensure accurate convergence to the correct solution.

Finally, Fig. \ref{fig:Null-space-BE} (c) showcases the 
immediate increase and subsequent slow decay in the solution's rank as the simulation progresses. This increase demonstrates the algorithm's effectiveness in augmenting the solution rank until it can accommodate all the dominant modes supported by the implicit integrator for a given time-step size. 

\subsection{Two-species LBFP Thermal Equilibration problem}
We test next our algorithm on the nonlinear LBFP model presented in section \ref{sec:fp}. Similarly to the heat equation results, we will show that our algorithm yields a temporal error comparable to that of the full-rank integrator. Importantly, our model succeeds in maintaining a low-rank structure consistently throughout the simulation, resulting in tremendous computational speedup. 

In our study, we non-dimensionalize the initial parameters using reference values detailed in \cite{taitano2015charge}. Our simulations explore the collision dynamics between ions and electrons with a realistic mass ratio, $m_i=1$ for ions and $m_e=\frac{1}{1836}$ for electrons. 
The initial particle distributions for each species are described by the following bi-Maxwellian initial condition:
\begin{equation}
f(v_{\alpha,1}, v_{\alpha,2}, t=0) = \frac{n_{\alpha,0}}{{2\pi v_{th,\alpha} ^2} } \left(0.5e^{-\frac{m_\alpha ( (v_{\alpha,1} - {u}^1_{\alpha,0})^2 + (v_{\alpha,2} - {u}^2_{\alpha,0})^2 )}{2 T_{\alpha,0}}} + 0.5e^{-\frac{m_\alpha ( (v_{\alpha,1} + {u}^1_{\alpha,0})^2 + (v_{\alpha,2} + {u}^2_{\alpha,0})^2 )}{2 T_{\alpha,0}}}\right)
\label{eq:initial_cond}
\end{equation}
In this model, the terms $n_{\alpha,0}$, ${u}^1_{\alpha,0}$, ${{u}}^2_{\alpha,0}$, and $v_{th_{\alpha,0}}=\sqrt{T_{\alpha,0}/{m_\alpha}}$ represent the initial density, the initial drift velocities in the $v_{\alpha,1}$ and $v_{\alpha,2}$ directions, and the initial thermal velocity of species $\alpha$, respectively, at time $t=0$. We define the initial drift velocities as ${u}^1_{i,0} = {u}^2_{i,0} = 2$ for ions and ${u}^1_{e,0} = {u}^2_{e,0} = 10$ for electrons. The initial temperatures are set to $T_{i,0}=1.1$ for ions and $T_{e,0}=0.9$ for electrons. The initial thermal velocities ($v_{th,i}$ and $v_{th,e}$) are $1.048$ for ions and $40.6497$ for electrons. The initial total drift velocities 
are zero for both species. The initial total ion temperature (i.e., including both Maxwellians) is \(5.1\), whereas the initial total electron temperature is \(0.954466\). 

We factorize the initial distribution function as $U_\alpha S_\alpha V_\alpha^T = \texttt{svd}\left(f(v_{\alpha,1}, v_{\alpha,2}, t=0)\right)$. Since the chosen initial condition is of rank two, we retain the two foremost singular values contained within $S_\alpha(1:2,1:2)$, along with their corresponding singular vectors $U_\alpha(:,1:2)$ and $V_\alpha(:,1:2)$. We employ an SVD truncation threshold of $10^{-8}\sigma_1$, where $\sigma_1$ is the largest singular value of the solution. This threshold remains fixed unless otherwise stated.
This approach truncates the initial condition to a considerably reduced dimensional space compared to the full tensor grid, thus starting the simulation with a compact and computationally efficient representation. 
We set the simulation domains to span $[-10v_{th,i}, 10v_{th,i}]$ in both directions for ions, and similarly $[-10v_{th,e}, 10v_{th,e}]$ for electrons. 
As simulations progress, the system temperatures converge to an equilibrium temperature, $\Bar{T}= 3.02723$, as described by the analytical expression, Eq. \eqref{eq:Temp-eq}. For each velocity direction, we utilize a nominal grid resolution of $N_{v,\alpha}=1000$, leading to a total of $10^6$ unknowns in our solution matrix for the full-rank representation. The residual tolerance for our simulations 
is set to be proportional to the LTE of the utilized DIRK scheme, i.e. $\epsilon^{(k)}_{tol}=\text{C}_k \Delta t^{(k+1)}$, where $p$ is the order of the temporal integration utilized and $C_k$ are user-specified constants. For BE simulations, $\text{C}_1=1$; for DIRK2, $\text{C}_1=\text{C}_2=10^{-3}$; and for DIRK3 simulations, $\text{C}_1=\text{C}_2=\text{C}_3=10^{-3}$. The time step $\Delta t$ is proportional to the square of the velocity grid size for ion, denoted as $\Delta v_i$, with 
$\lambda=\frac{\Delta t}{{\Delta v_i} ^2}$ ranging from 100 to 900. 

\begin{figure}
    \centering
\includegraphics[width=\textwidth]{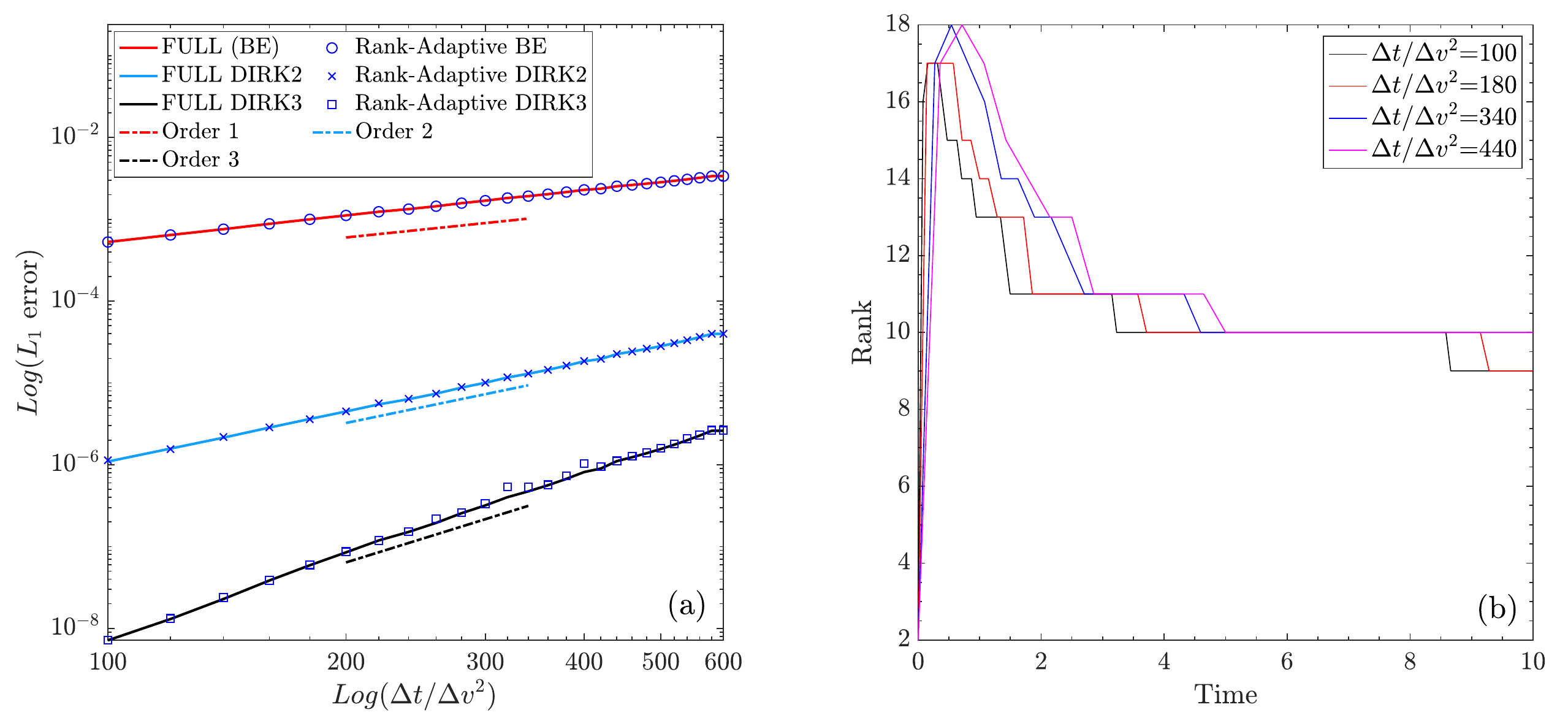}
    \caption{Simulation of the Fokker-Planck equation employing three temporal integrators--specifically Backward Euler, DIRK2, and DIRK3. In Fig. (a), we illustrate a temporal error convergence study spanning a range $\lambda=\frac{\Delta t}{{\Delta v}^2}$ values [100:600] for both the low-rank integrator and the classical full-rank counterpart. Fig. (b) shows the rank evolution as a function of time for the Backward Euler simulation.}
    \label{FP_fig1}
\end{figure}

In Fig. \ref{FP_fig1} (a), we present results from simulations conducted using Algorithm \ref{DIRK-LR} for the Butcher table corresponding to three different integrators: Backward Euler, DIRK 2, and DIRK 3. Our proposed algorithm demonstrates error convergence identical to that of the full-rank integrator while achieving remarkable speedup, discussed below. Additionally, the algorithm adeptly captures the rank evolution of the solution, as illustrated in Fig. \ref{FP_fig1} (b). 

\begin{figure}
    \centering
\includegraphics[width=\textwidth]{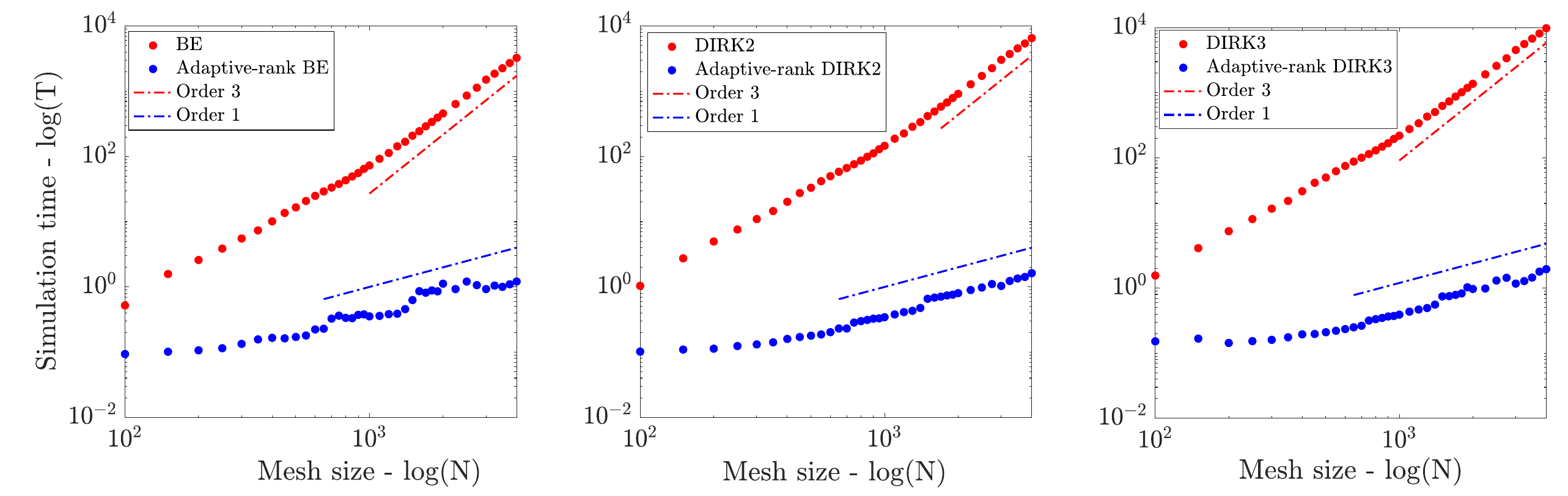}
    \caption{ Comparison of computational complexity between the adaptive-rank integrator (illustrated in blue) and the full-rank integrator (shown in red) for BE, DIRK2, and DIRK3. Here, $N$ represents the number of grid points for each velocity dimension, denoted as $N_{v_1}$ and $N_{v_2}$, where for simplicity we set $N=N_{v_1}=N_{v_2}$. The simulation time was measured using MATLAB's \texttt{timeit} function.}
    \label{fig: complexity_fp}
\end{figure}

As discussed in Section \ref{sec: 2.3}, our algorithm scales linearly with the 1D mesh dimension, i.e., $O(N)$, in contrast to the $O(N^3)$ computational scaling of the full-rank integrator. This significant reduction in computational complexity underscores the potential efficiency gains achievable through exploiting the low-rank structure of the solution. We adopt the same setup as before, except that we fix the time-step at $\Delta t =0.1$,
we vary the mesh resolutions from $100-4000$.  Fig. \ref{fig: complexity_fp} demonstrates the $\mathcal{O}(N)$ computational scaling of our low-rank algorithm and the O($N^3$) of the full-rank classical Sylvester solver with the same integrators. 

\begin{figure}
    \centering
\includegraphics[width=0.7\textwidth]{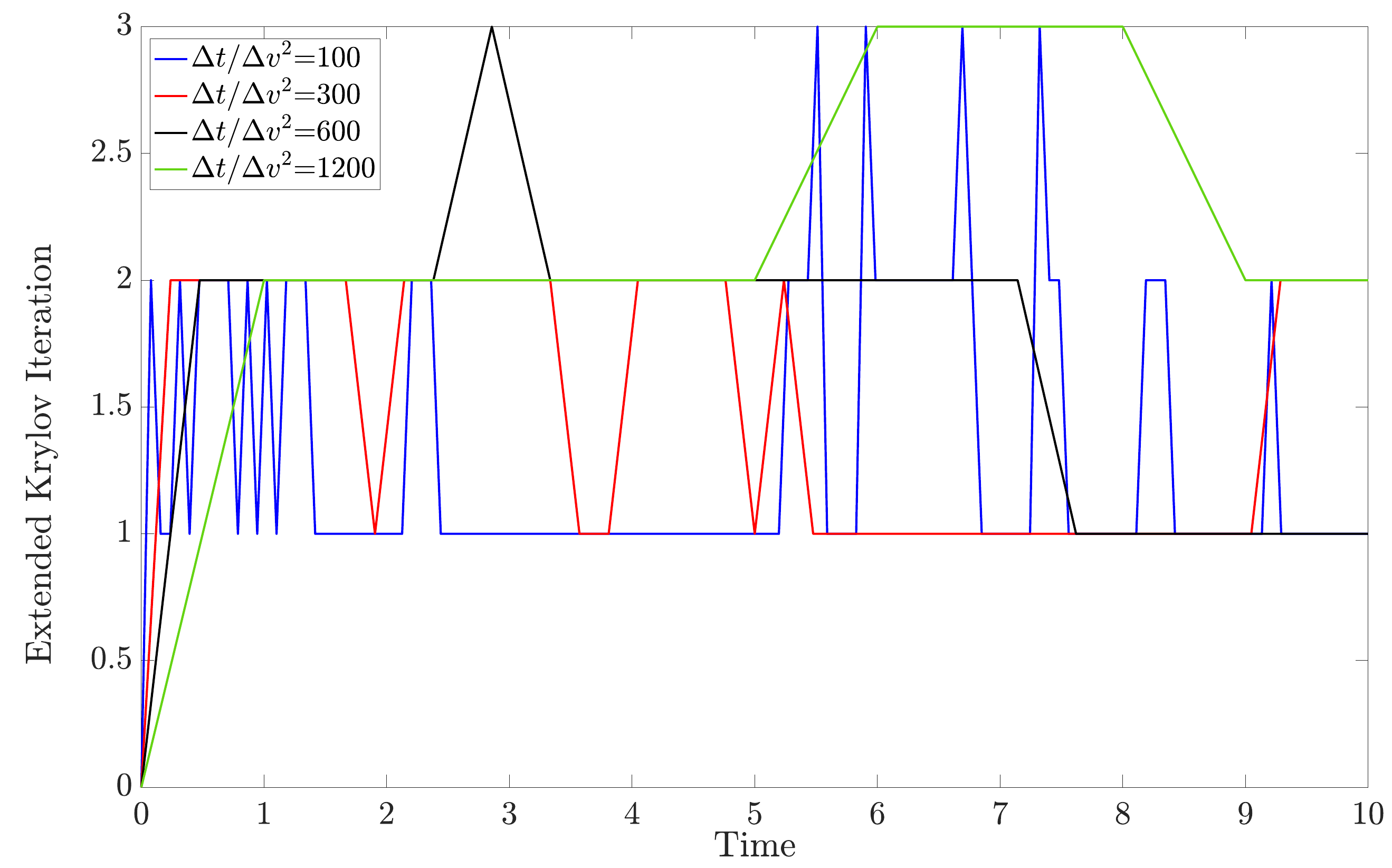}
    \caption{Evolution of the extended Krylov iterations  per timestep over time for various $\lambda$ numbers for BE. This iteration number is representative for each stage in DIRK.}
    \label{fig:comp_FP_LR}
\end{figure}

Next, to numerically show that the rank $r_m$ 
remains low across the simulations,  we record the number of Krylov iterations, $m$, for a range of $\lambda$ values. Fig. \ref{fig:comp_FP_LR} shows that this quantity remains bounded even for relatively larger time steps. This attests to the fast convergence property of the extended Krylov subspace. 

\begin{figure}
    \centering
\includegraphics[width=\textwidth]{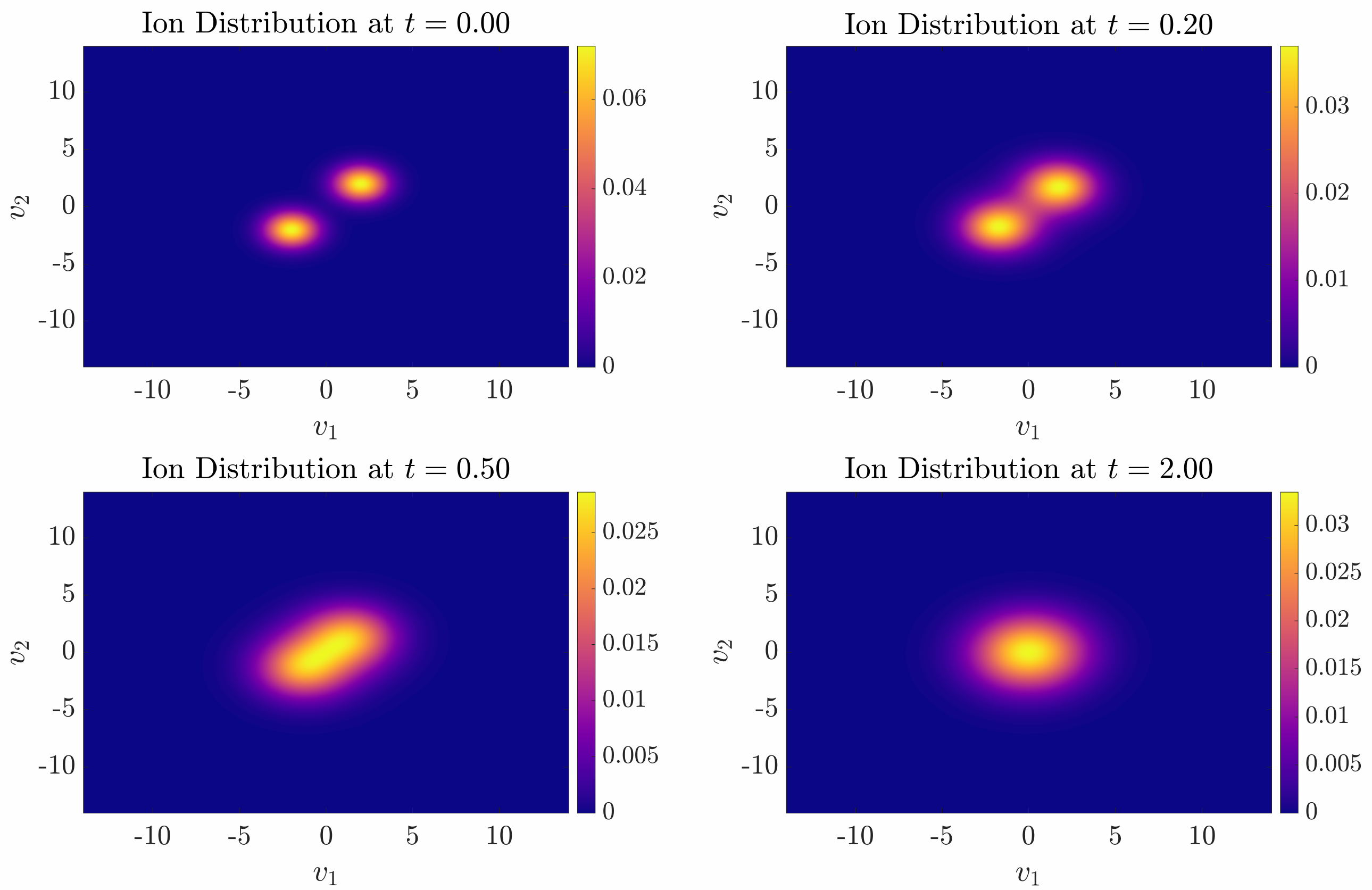}
    \caption{Temporal evolution of ion distribution 
    at selected time intervals: $t=0$, $t=0.2$, $t=0.5$, and $t=2$. The distribution converges to a single Maxwellian, indicating the system's transition to an equilibrium.}
    \label{fig:LBFP_ions}
\end{figure}

\begin{figure}
    \centering
\includegraphics[width=0.9\textwidth]{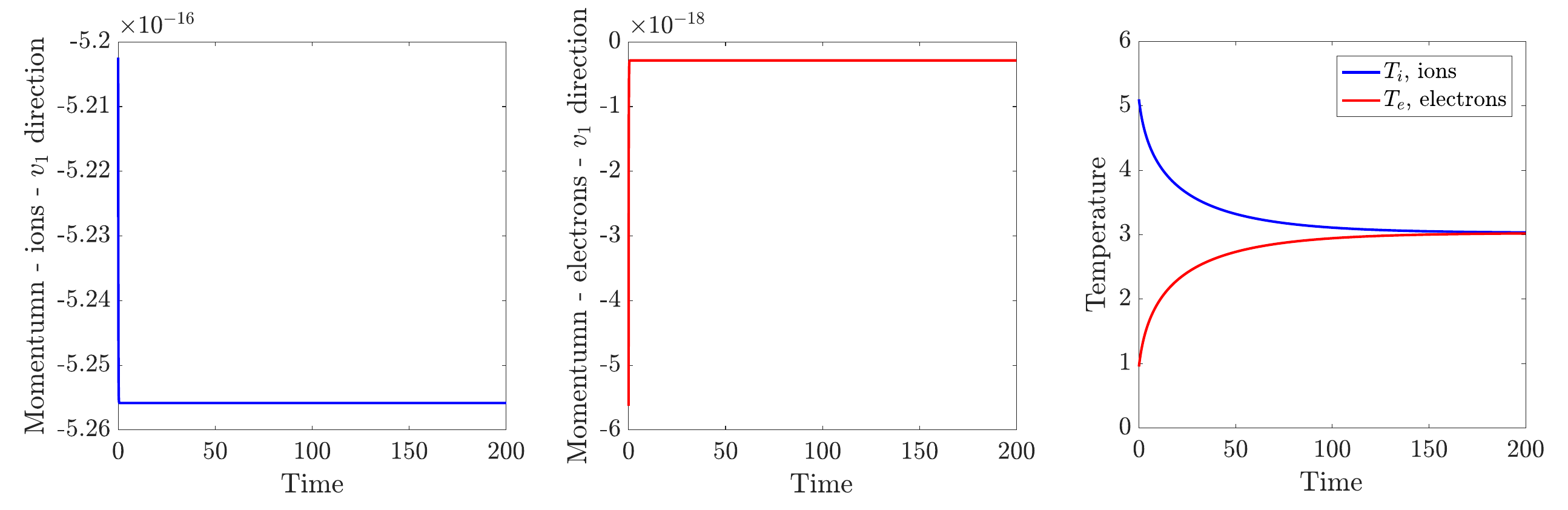}
    \caption{Temporal history of $v_1$-component of total momentum and temperature, for both ions and electrons, for the two-species relaxation test.}
    \label{fig:ODE_evolution}
\end{figure}

\begin{figure}
    \centering
\includegraphics[width=0.9\textwidth]{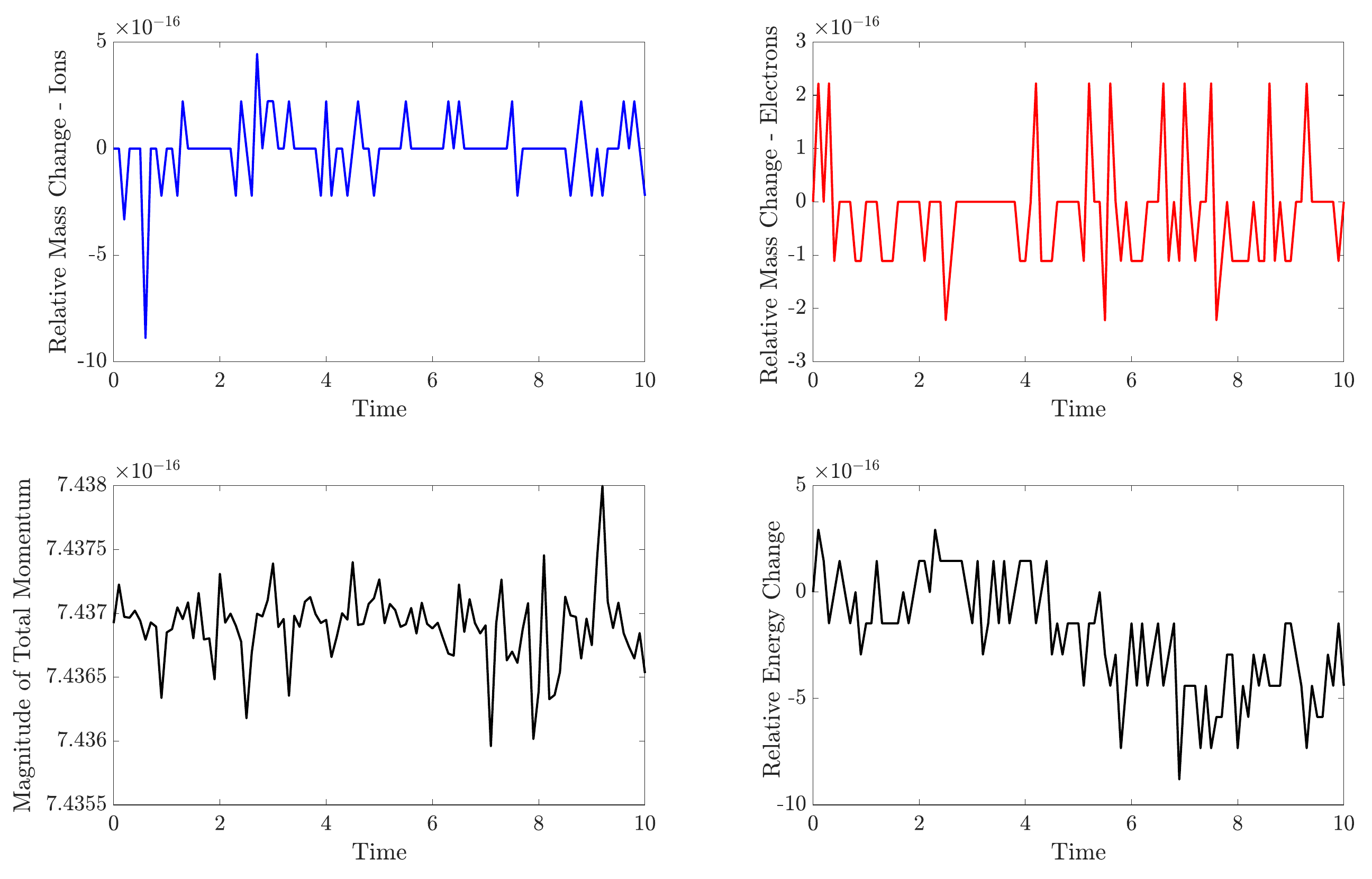}
    \caption{Temporal history of conservation errors from mass (for both ions and electrons), magnitude of total momentum, and total energy for the two-species relaxation test.}
    \label{fig:conservation}
\end{figure}

Finally, we present the solution snapshots at different times and evaluate the macroscopic conservation properties of the proposed algorithm with LoMaC procedure. We integrate the solution to a final time $T_f=10$ with $\Delta t =0.1$. Fig. \ref{fig:LBFP_ions} presents the time evolution of ion distributions at four distinct snapshots: $t=0$, $t=0.2$, $t=0.5$, and $t=2$. Initially, the ion distribution is characterized by two Maxwellian profiles, centered at coordinates $(-2, -2)$ and $(2, 2)$. As time progresses, the Maxwellians undergo transport and diffusion, eventually settling into a single, stable Maxwellian distribution. Fig. \ref{fig:ODE_evolution} illustrates the time history of mass, momentum, and energy of each species computed per Eqs. \eqref{continuityODE}-\eqref{energyODE}. At the beginning of the simulation, we observe a rapid increase in electron temperature due to the conversion of kinetic energy into internal energy, leading to an early equilibration of momentum. Following this momentum equilibration, the slower process of temperature relaxation dominates, resulting in a gradual decrease in ion temperature and a corresponding increase in electron temperature until a steady-state temperature equilibrium is achieved.
Fig. \ref{fig:conservation} illustrates the evolution of the relative mass and energy variations from that of the initial conditions, alongside the norm of the error of the momentum vector. These results confirm that the proposed algorithm conserves the first three moments of the LBFP equation with machine precision. 

\section{Conclusions}
\label{sec:conclusions}

We have proposed a high-order adaptive-rank implicit integrator that leverages the power of extended Krylov-subspace methods for dynamic rank discovery. Our emphasis has been on high-order temporal accuracy (via DIRK schemes), strict conservation properties (via the LoMaC projection), and superoptimal computational complexity (i.e., that does not scale as the total number of unknowns, $N^d$, with $d$ the number of dimensions). We have applied our algorithm to the classical heat equation and the more challenging Lenard-Bernstein Fokker-Planck (LBFP) nonlinear equation. The scheme has demonstrated 1) the ability to discover rank evolution faithfully and efficiently, 2) high-order temporal accuracy (up to third order in this study), 3) strict conservation properties (to numerical round-off), and 4) superoptimal computational complexity scaling as $d N$, thereby breaking the curse of dimensionality. Future extensions of this study include extension to higher dimensions (which will require tensor factorization techniques), and applications to more general nonlinear kinetic models. 

\section*{Acknowledgments}

W.T. and L.C. and H.E. (during a summer internship in 2023) were supported by Triad National Security, LLC under contract 89233218CNA000001 and DOE Office of Applied Scientific Computing Research (ASCR) through the Mathematical Multifaceted Integrated Capability Centers program. J.Q. and H.E. were supported by NSF grant NSF-DMS-2111253, J.Q. was supported by Air Force Office of Scientific Research FA9550-22-1-0390 and Department of Energy DE-SC0023164.

\bibliographystyle{abbrv}
\bibliography{sample}

\appendix
\section{Chang-Cooper equilibrium preserving discretization of the LBFP collision operator}
\label{app:chang-cooper}
We describe the construction of the differentiation matrices \(A_1\) and \(A_2\) for the LBFP introduced in step 2, Section 3.2. We denote the discrete solution of the equation \eqref{eq:lbfp} at specific grid points by \(f(v_1(i), v_2(j)) = \bF_{i,j}\), and 
the coordinates \(\{v_1(i), v_2(j)\}\) by \(\{v_{i}^{1}, v^{2}_{j}\}\) for conciseness. For the sake of simplicity, we assume \(v_1 = v_2\), implying that the velocity grid in the first dimension is identical to that of the second dimension. Additionally, we assume $\Delta v_1 = \Delta v_2=\Delta v$, indicating the grid spacing in both dimensions is identical. These assumptions simplify our discussion without loss of generality.

We apply the MOL for discretizing equation  \eqref{eq:lbfp}. The MOL approach considers the discretization of the velocity variables, leading to a system of ODEs with respect to time. 
Let's define the collisional flux vector at the faces of the cell $(i,j)$ as:
\begin{align}
    {{\vec{\phi}}_{i,j}}=[\phi_{i+\frac{1}{2},j}^1,\phi_{i,j+\frac{1}{2}}^2]= \biggl[ 
    &\sum^{N_s}_{{\beta=1}} \nu_{\alpha\beta} \left(    D_{\alpha\beta}\frac{F_{i+1,j} - F_{i,j} }{\Delta v}  + \left(v_{i+\frac{1}{2}} - {u^1}_{\alpha\beta} \right){F_{i+\frac{1}{2},j}} \right),\nonumber\\
&\sum^{N_s}_{{\beta=1}} \nu_{\alpha\beta} 
\left(D_{\alpha\beta}\frac{F_{i,j+1} - F_{i,j} }{\Delta v} + \left(v_{j+\frac{1}{2}} - {u^2}_{\alpha\beta} \right){F_{i,j+\frac{1}{2}}}\right) 
    \biggl]
\end{align}
The discretization of the LBFP operator for species $\alpha$ at the cell $(i,j)$ then reads: 
\begin{equation}
\frac{\partial F_{i,j}}{\partial t} = \frac{\phi^1_{i+1/2,j}-\phi^1_{i-1/2,j}}{\Delta v} + \frac{\phi^2_{i,j+1/2}-\phi^2_{i,j-1/2}}{\Delta v}.\label{eq:flux_form}
\end{equation}
We impose zero-flux at the boundaries of the domains, i.e. $\phi^1_{\frac{1}{2},j}=0$, $\phi^1_{N_v+\frac{1}{2},j}=0$,  $\phi^2_{i,\frac{1}{2}}=0$, and $\phi^2_{i,N_v+\frac{1}{2}}=0$. 
We use a weighted average of the neighboring known cell values to approximate the unknown half-cell values $F_{i,j+\frac{1}{2}}$, $F_{i,j-\frac{1}{2}}$, $F_{i+\frac{1}{2},j}$, and $F_{i-\frac{1}{2},j}$: 
\begin{align*}
F_{i,j+\frac{1}{2}} &= \delta_j F_{i,j} + (1-\delta_j)F_{i,j+1}, \\
F_{i+\frac{1}{2},j} &= \delta_i F_{i,j} + (1-\delta_i)F_{i+1,j}.
\end{align*}
Here, the weights $\delta_j$ and $\delta_i$ are constrained such that $0 \leq \{\delta_j,\delta_i\} \leq \frac{1}{2}$. Typically, a standard second-order finite-difference scheme would set $\delta_j = \frac{1}{2}$ and $\delta_i = \frac{1}{2}$ to approximate these quantities.

However, in their seminal work, Chang and Cooper \cite{chang1970practical} identified that this conventional choice leads to a scheme that does not preserve the Maxwellian equilibrium. To address this limitation, they proposed the following alternative formulation for $\delta_{\alpha\beta,m}$:
\[
\delta^{k}_{\alpha\beta,m} = \frac{1}{w^{k}_{\alpha\beta,m}} - \frac{1}{\exp(w^{k}_{\alpha\beta,m}) - 1},
\]
where the superscript $k$ denotes the velocity direction and underscript $m$ the discretization node, and $w^{k}_{\alpha\beta,m}$ is defined as
\[
w^{k}_{\alpha\beta,m} = \frac{\Delta v(v_{m+\frac{1}{2}} - u^{k}_{\alpha\beta})}{D_{\alpha\beta}}.
\]
Note the addition of a subscript $\alpha\beta$, since the weighting parameter in our expression would depend  on the inter-species drift velocities. This choice of $\delta^k_{\alpha\beta,m}$ is designed to ensure the preservation of the Maxwellian equilibrium distribution.


Expanding the divergence operator in equation \eqref{eq:flux_form}, we have:
\begin{align*}
    \frac{\partial F_{i,j}}{\partial t} &= \sum_{\beta=1}^{N_s} \gamma_{\alpha\beta} \Bigg( 
    D_{\alpha\beta} \frac{F_{i+1,j} - 2F_{i,j} + F_{i-1,j}}{\Delta {v}^2} +
    D_{\alpha\beta} \frac{F_{i,j+1} - 2F_{i,j} + F_{i,j-1}}{\Delta {v}^2} \\
    &\quad +  \frac{\left( v_{i+\frac{1}{2}} - u^{x}_{\alpha,\beta} \right)F_{i+\frac{1}{2},j} - \left( v_{i-\frac{1}{2}} - u^{x}_{\alpha,\beta} \right)F_{i-\frac{1}{2},j}}{\Delta v} +
     \frac{ \left( v_{j+\frac{1}{2}} - u^{y}_{\alpha,\beta} \right)F_{i,j+\frac{1}{2}} - \left( v_{j-\frac{1}{2}} - u^{y}_{\alpha,\beta} \right)F_{i,j-\frac{1}{2}}}{\Delta v} \Bigg),
\end{align*}
Substituting the cell averages, boundary conditions, and collecting terms, we find:

\begin{align}
    \frac{\partial F_{i,j}}{\partial t} = & F_{i+1,j} \underbrace {\left( \sum_{\beta=1}^{N_s} \gamma_{\alpha\beta} \left( \frac{D_{\alpha\beta}}{{\Delta v}^2} + \frac{1}{\Delta v}(1 - \delta^{1}_{\alpha\beta,i}) (v_{i+\frac{1}{2}} - u_{1,\alpha\beta}) \right) \right) }_{s_{i+1}}   \label{App:LBFP}\\
    + & F_{i,j}  \underbrace {\left( \sum_{\beta=1}^{N_s} \gamma_{\alpha\beta} \left( \frac{-2D_{\alpha\beta}}{{\Delta v}^2} + \frac{1}{\Delta v} \delta^{1}_{\alpha\beta,i} (v_{i+\frac{1}{2}} - u_{1,\alpha\beta}) + \frac{1}{\Delta v} (1 - \delta^{1}_{\alpha\beta,i-1}) (v_{i-\frac{1}{2}} - u_{1,\alpha\beta}) \right) \right)}_{d_i}\nonumber \\
    + & F_{i-1,j}  \underbrace {\left( \sum_{\beta=1}^{N_s} \gamma_{\alpha\beta} \left( \frac{D_{\alpha\beta}}{{\Delta v}^2} - \frac{1}{\Delta v} \delta^{1}_{\alpha\beta,i-1} (v_{i-\frac{1}{2}} - u_{1,\alpha\beta}) \right) \right) }_{l_{i-1}}\nonumber\\
    + & F_{i,j+1}  \underbrace {\left( \sum_{\beta=1}^{N_s} \gamma_{\alpha\beta} \left( \frac{D_{\alpha\beta}}{{\Delta v}^2} + \frac{1}{\Delta v} (1 - \delta^{2}_{\alpha\beta,j}) (v_{j+\frac{1}{2}} - u_{1,\alpha\beta}) \right) \right)}_{\Tilde{s}_j} \nonumber\\
    + & \underbrace {F_{i,j} \left( \sum_{\beta=1}^{N_s} \gamma_{\alpha\beta} \left( \frac{-2D_{\alpha\beta}}{{\Delta v}^2} + \frac{1}{\Delta v} \delta^{2}_{\alpha\beta,j} (v_{j+\frac{1}{2}} - u_{1,\alpha\beta}) + \frac{1}{\Delta v} (1 - \delta_{2,j-1}) (v_{j-\frac{1}{2}} - u_{1,\alpha\beta}) \right) \right) }_{\Tilde{d}_j}\nonumber\\
    + & F_{i,j-1} \underbrace {\left( \sum_{\beta=1}^{N_s} \gamma_{\alpha\beta} \left( \frac{D_{\alpha\beta}}{{\Delta v}^2} - \frac{1}{\Delta v} \delta^{2}_{\alpha\beta,j-1} (v_{j-\frac{1}{2}} - u_{1,\alpha\beta}) \right) \right)}_{\Tilde{l}_{j-1}},\nonumber
\end{align}
which provides the definitions of the discretization matrix elements.
The zero-flux boundary conditions yield the following relation for the ghost cells $F_{0,j}, F_{N+1,j}$, $F_{i,0}$, and $F_{i,N_v+1}$:
\begin{equation}
    F_{0,j}=\frac{ s_1}{ l_0}F_{1,j}, \quad      F_{N_v+1,j}=\frac{ l_{N_v}}{ s_{N_v+1}}F_{N_v,j}, \quad F_{i,0}=\frac{ \Tilde{s}_1}{ \Tilde{l}_0}F_{i,1}, \quad      F_{i,N_v+1}=\frac{ \Tilde{l}_{N_v}}{ \Tilde{s}_{N_v+1}}F_{i,N_v}.
\end{equation}

The resulting matrix differential equation from the MOL above reads:
\begin{equation}
\frac{\partial \bF}{\partial t}= A_1 \bF + \bF {A_2}^T
\end{equation}
where \(A_1\) and \(A_2\) are tridiagonal matrices defined as follows:
\[
A_1 = \begin{bmatrix}
d_{1}+s_{1} & s_2 & 0 & \cdots & 0 \\
l_1 & d_2 & s_3 & \ddots & \vdots \\
0 & \ddots & \ddots & \ddots & 0 \\
\vdots & \ddots & l_{N_v -1} & d_{N_v -1} & s_{N_v} \\
0 & \cdots & 0 & l_{N_v -1}& d_{N_v}+ l_{N_v} \\
\end{bmatrix},
\]
with \(d_i\), \(s_i\), and \(l_i\) derived from the coefficients multiplying \(F_{i,j}\), \(F_{i+1,j}\), and \(F_{i-1,j}\) respectively.
\[
A_2 = \begin{bmatrix}
\Tilde{d}_{1}+ \Tilde{s_1} & \Tilde{s}_2 & 0 & \cdots & 0 \\
\Tilde{l}_1 & \Tilde{d}_3 & \Tilde{s}_2 & \ddots & \vdots \\
0 & \ddots & \ddots & \ddots & 0 \\
\vdots & \ddots & \Tilde{l}_{N_v} & \Tilde{d}_{N_v -1} & \Tilde{s}_{N_v} \\
0 & \cdots & 0 & \Tilde{l}_{N_v-1} & \Tilde{d}_{N_v} + \Tilde{l}_{N_v}  \\
\end{bmatrix},
\]
with \(\Tilde{d}_i\), \(\Tilde{s}_i\), and \(\Tilde{l}_i\) derived from the coefficients affecting \(F_{i,j}\), \(F_{i,j+1}\), and \(F_{i,j-1}\) respectively. 
Note that the elements of these matrices depend on the moment quantities per the Chang-Cooper interpolation. In the context of the DIRK method, where several time stages are considered, it raises the question as to how to incorporate these moments in time for the velocity discretization. One option is to use the moments at the same temporal location per stage. However, we have found this leads to temporal order reduction. 
Another option, which we employ here, is to use the moments at the final DIRK stage (which we know a priori since we integrate the moment ODEs first) for all stages. Since these moments are used for equilibrium preservation of the velocity-space scheme and do not determine temporal accuracy, we have not seen temporal order reduction 
(as evidenced by the numerical results).

\end{document}